\documentclass{article}%
\usepackage{amsmath}
\usepackage{amsfonts}
\usepackage{amssymb}
\usepackage{graphicx}%
\providecommand{\U}[1]{\protect\rule{.1in}{.1in}}

\newenvironment{proof}[1][Proof]{\noindent\textbf{#1.} }{\ \rule{0.5em}{0.5em}}
\begin{document}

\title{Asymptotic behaviour of a differential operator with a finite number of
transmission conditions}
\author{{\Large Erdo\u{g}an \c{S}en and O. Sh. Mukhtarov$^*$}}
\date{}
\maketitle

{\scriptsize {Department of Mathematics Engineering, Istanbul Technical
University, Maslak, 34469 Istanbul, Turkey}}

{\scriptsize {Department of Mathematics , Istanbul Technical
University, Maslak, 34469 Istanbul, Turkey}}

{\scriptsize {$^*$Department of Mathematics, Faculty of Arts and Science, Gaziosmanpa\c{s}a University,\\
 60250 Tokat, Turkey.}}

{\scriptsize e-mail: erdogan.math@gmail.com,
$^*$omukhtarov@yahoo.com}



ABSTRACT

In this paper following the same methods in [M. Kadakal, O. Sh. Mukhtarov,
Sturm-Liouville problems with discontinuities at two points, Comput. Math.
Appl., 54 (2007) 1367-1379] we investigate discontinuous two-point boundary
value problems with eigenparameter in the boundary conditions and with
transmission conditions at the finitely many points of discontinuity. A
self-adjoint linear operator $A$ is defined in a suitable Hilbert space $H$
such that the eigenvalues of such a problem coincide with those of $A$. We
obtain asymptotic formulas for the eigenvalues and eigenfunctions. Also we
show that the eigenfunctions of $A$ are complete in $H$.

\bigskip

\section{Introduction}

The theory of discontinuous Sturm-Liouville type problems mainly has been
developed by Mukhtarov and his students (see [1-11]).It is well-known that
many topics in mathematical physics require the investigation of eigenvalues
and eigenfunctions of Sturm-Liouville type boundary value problems. In recent
years, more and more researchers are interested in the discontinuous
Sturm-Liouville problems and its applications in physics (see $[1-28]$).

Discontinuous Sturm-Liouville problems with supplementary transmission
conditions at the point(s) of discontinuity have been investigated in [6-9, 23-27].

In this study, we examine eigenvalues and eigenfunctions of the differential
equation%
\begin{equation}
\tau u:=-u^{\prime\prime}(x)+q(x)u(x)=\lambda u(x) \tag*{(1)}%
\end{equation}
on $\left[  -1,h_{1}\right)  \cup\left(  h_{1},h_{2}\right)  \cup
...\cup\left(  h_{m},1\right]  ,$ with boundary conditions%
\begin{align}
\tau_{1}u  &  :=\alpha_{1}u\left(  -1\right)  +\alpha_{2}u^{\prime
}(-1)=0,\tag*{(2)}\\
\tau_{2}u  &  :=\left(  \beta_{1}^{\prime}\lambda+\beta_{1}\right)  u\left(
1\right)  +\left(  \beta_{2}^{\prime}\lambda+\beta_{2}\right)  u^{\prime}(1)=0
\tag{3}%
\end{align}
and transmission conditions at the points of discontinuity $x=h_{i}$
$(i=\overline{1,m}),$%

\begin{align}
\tau_{2i+1}u  &  :=u\left(  h_{i}-0\right)  -\delta_{i}u\left(  h_{i}%
+0\right)  =0,\tag*{(4)}\\
\tau_{2i+2}u  &  :=u^{\prime}\left(  h_{i}-0\right)  -\delta_{i}u^{\prime
}\left(  h_{i}+0\right)  =0, \tag{5}%
\end{align}
where $-1<h_{1}<h_{2}<...<h_{m}<1,$ $q(x)$ is a given real-valued function
continuous in $\left[  -1,h_{1}\right)  ,\left(  h_{1},h_{2}\right)  ,...,$
$\left(  h_{m},1\right]  ~$and has finite limits $q(h_{i}\pm0)=\lim
_{x\rightarrow h_{i}\pm0}q(x)$ $\left(  i=\overline{1,m}\right)  $; $\lambda$
is a complex eigenvalue parameter; $\delta_{i}$ $\left(  i=\overline
{1,m}\right)  ,\alpha_{j},\alpha_{j}^{^{\prime}},\beta_{j},\beta_{j}%
^{^{\prime}}$ $\left(  j=1,2\right)  $ are real numbers; $\left\vert
\alpha_{1}\right\vert +\left\vert \alpha_{2}\right\vert \neq0$ and $\delta
_{i}\neq0$ $\left(  i=\overline{1,m}\right)  .$ As following $\left[
6\right]  $ every where we assume that $\rho=\beta_{1}^{^{\prime}}\beta
_{2}-\beta_{1}\beta_{2}^{^{\prime}}>0.$

\section{Operator Formulation}
By using the method, introduced in [6]  we shall define direct sum
of Hilbert spaces but with the usual inner product replaced by
appropriate multiples as follows.
In the Hilbert space $H=L_{2}\left(  -1,1\right)  \oplus%
\mathbb{C}
$ we define an inner product by%

\begin{align*}
\left\langle F,G\right\rangle  &  :=%
{\displaystyle\int\limits_{-1}^{h_{1}}}
f\left(  x\right)  \overline{g\left(  x\right)  }dx+\delta_{1}^{2}%
{\displaystyle\int\limits_{h_{1}}^{h_{2}}}
f\left(  x\right)  \overline{g\left(  x\right)  }dx+\delta_{1}^{2}\delta
_{2}^{2}%
{\displaystyle\int\limits_{h_{2}}^{h_{3}}}
f\left(  x\right)  \overline{g\left(  x\right)  }dx+\\
&  ...+%
{\displaystyle\prod\limits_{i=1}^{m}}
\delta_{i}^{2}%
{\displaystyle\int\limits_{h_{m}}^{1}}
f\left(  x\right)  \overline{g\left(  x\right)  }dx+\frac{%
{\displaystyle\prod\limits_{i=1}^{m}}
\delta_{i}^{2}}{\rho}f_{1}\overline{g_{1}}%
\end{align*}
for%
\[
F:=\left(
\begin{array}
[c]{c}%
f\left(  x\right) \\
f_{1}%
\end{array}
\right)  ,\text{ }G:=\left(
\begin{array}
[c]{c}%
g\left(  x\right) \\
g_{1}%
\end{array}
\right)  \in H
\]
following $\left(  7\right)  $ for convenience we put%
\begin{align*}
R_{1}\left(  u\right)   &  :=\beta_{1}u\left(  1\right)  -\beta_{2}%
u^{^{\prime}}\left(  1\right)  ,\\
R_{1}^{^{\prime}}\left(  u\right)   &  :=\beta_{1}^{^{\prime}}u\left(
1\right)  -\beta_{2}^{^{\prime}}u^{^{\prime}}\left(  1\right)  .
\end{align*}
For functionals $f\left(  x\right)  ,$ which defined on $\left[
-1,h_{1}\right)  \cup\left(  h_{1},h_{2}\right)  \cup...\cup\left(
h_{m},1\right]  $ and has finite limits $f\left(  h_{i}\pm0\right)
:=\lim_{x\rightarrow h_{i}\pm0}f\left(  x\right)  $ $\left(  i=\overline
{1,m}\right)  $. By $f_{i}\left(  x\right)  $ $\left(  i=\overline
{1,m+1}\right)  $ we denote the functions%
\begin{align*}
f_{1}\left(  x\right)   &  :=\left\{
\begin{array}
[c]{ll}%
f(x), & x\in\lbrack-1,{h}_{1}),\\
\lim_{x\rightarrow h_{1}-0}f(x), & x=h_{1},
\end{array}
\right.  \text{ \ }f_{2}\left(  x\right)  :=\left\{
\begin{array}
[c]{cc}%
\lim_{x\rightarrow h_{1}+0}f(x), & x=h_{1},\\
f(x), & x\in\left(  h_{1},{h}_{2}\right)  ,\\
\lim_{x\rightarrow h_{2}-0}f(x), & x=h_{2},
\end{array}
...\right. \\
,f_{m}\left(  x\right)   &  :=\left\{
\begin{array}
[c]{cc}%
\lim_{x\rightarrow h_{m-1}+0}f(x), & x=h_{m-1},\\
f(x), & x\in\left(  h_{m-1},{h}_{m}\right)  ,\\
\lim_{x\rightarrow h_{m}-0}f(x), & x=h_{m},
\end{array}
\right.  \text{ \ }f_{m+1}\left(  x\right)  :=\left\{
\begin{array}
[c]{ll}%
\lim_{x\rightarrow h_{m}+0}f(x), & x=h_{m},\\
f(x), & x\in\left(  h_{m},{1}\right]
\end{array}
\right.
\end{align*}
which are defined on $\Omega_{1}:=\left[  -1,h_{1}\right]  ,$ $\Omega
_{2}:=\left[  h_{1},h_{2}\right]  ,...,\Omega_{m}:=\left[  h_{m-1}%
,h_{m}\right]  ,$ $\Omega_{m+1}:=\left[  h_{m},1\right]  $ respectively. We
can rewrite the considered problem $\left(  1\right)  -\left(  5\right)  $ in
the operator formulation as%
\[
AF=\lambda F
\]
where%
\[
F:=\left(
\begin{array}
[c]{c}%
f\left(  x\right) \\
-R_{1}^{^{\prime}}\left(  f\right)
\end{array}
\right)  \in D\left(  A\right)
\]
and with%
\[%
\begin{array}
[c]{c}%
D\left(  A\right)  :=\left\{  F\in H\left\vert f_{i}\left(  x\right)  ,\text{
}f_{i}^{^{\prime}}\left(  x\right)  \text{ are absolutely continuous in
}\Omega_{i}\text{ }\left(  i=\overline{1,m+1}\right)  \text{, }\tau
f\in\right.  \right. \\
\left.
\begin{array}
[c]{c}%
L^{2}\left[  -1,1\right]  ,\text{ }\tau_{2i+1}u:=u\left(  h_{i}-0\right)
-\delta_{i}u\left(  h_{i}+0\right)  =0,\text{ }\tau_{2i+2}u:=u^{\prime}\left(
h_{i}-0\right)  -\delta_{i}u^{\prime}\left(  h_{i}+0\right)  =0\\
\left(  i=\overline{1,m+1}\right)  \text{ and }f_{1}=R_{1}^{^{\prime}}\left(
f\right)
\end{array}
\right\}
\end{array}
\]
and%
\[
AF=\left(
\begin{array}
[c]{c}%
\tau f\\
-R_{1}\left(  f\right)
\end{array}
\right)  .
\]
Consequently, the problem $(1)-(5)$ can be considered as the eigenvalue
problem for the operator $A$. Obviously, we have

\textbf{Lemma 2.1. }\textit{The eigenvalues of the boundary value problem
(1)-(5) coincide with those of }$A$\textit{, and its eigenfunctions are the
first components of the corresponding eigenfunctions of }$A.$

\textbf{Lemma 2.2.\ }\textit{The domain }$D\left(  A\right)  $\textit{ is
dense in }$H$.

\begin{proof}
Let $F=\left(
\begin{array}
[c]{c}%
f(x)\\
\text{ }f_{1}%
\end{array}
\right)  \in H$, $F\bot D(A)$ and $\widetilde{C}_{0}^{\infty}$ be a functional
set such that%
\[
\phi\left(  x\right)  =\left\{
\begin{array}
[c]{c}%
\phi_{1}(x),x\in\left[  -1,h_{1}\right)  ,\\
\phi_{2}(x),x\in\left(  h_{1},h_{2}\right)  ,\\
\vdots\\
\phi_{m+1}(x),x\in\left(  h_{m},1\right]
\end{array}
\right.
\]
for $\phi_{1}(x)\in\widetilde{C}_{0}^{\infty}\left[  -1,h_{1}\right)  $,
$\phi_{2}(x)\in\widetilde{C}_{0}^{\infty}\left(  h_{1},h_{2}\right)
$,...,$\phi_{m+1}(x)\in\widetilde{C}_{0}^{\infty}\left(  h_{m},1\right]  .$
Since $\widetilde{C}_{0}^{\infty}\oplus0\subset D(A)$ $\left(  0\in%
\mathbb{C}
\right)  $, any $U=\left(
\begin{array}
[c]{c}%
u(x)\\
0
\end{array}
\right)  \in\widetilde{C}_{0}^{\infty}\oplus0$ is orthogonal to $F$, namely%
\begin{align*}
\left\langle F,U\right\rangle  &  =%
{\displaystyle\int\limits_{-1}^{h_{1}}}
f\left(  x\right)  \overline{u\left(  x\right)  }dx+\delta_{1}^{2}%
{\displaystyle\int\limits_{h_{1}}^{h_{2}}}
f\left(  x\right)  \overline{u\left(  x\right)  }dx+\delta_{1}^{2}\delta
_{2}^{2}%
{\displaystyle\int\limits_{h_{2}}^{h_{3}}}
f\left(  x\right)  \overline{u\left(  x\right)  }dx+\\
&  ...+%
{\displaystyle\prod\limits_{i=1}^{m}}
\delta_{i}^{2}%
{\displaystyle\int\limits_{h_{m}}^{1}}
f\left(  x\right)  \overline{u\left(  x\right)  }dx=\left\langle
f,u\right\rangle _{1}.
\end{align*}
We can learn that $f(x)$ is orthogonal to $\widetilde{C}_{0}^{\infty}$ in
$L^{2}\left[  -1,1\right]  $, this implies $f(x)=0$. So for all $G=\left(
\begin{array}
[c]{c}%
g(x)\\
g_{1}%
\end{array}
\right)  \in D(A),$ $\left\langle F,G\right\rangle =\frac{%
{\displaystyle\prod\limits_{i=1}^{m}}
\delta_{i}^{2}}{\rho}f_{1}\overline{g_{1}}=0$. Thus $f_{1}=0$ since
$g_{1}=R_{1}^{\prime}\left(  g\right)  $ can be chosen arbitrarily. So
$F=\left(
\begin{array}
[c]{c}%
0\\
0
\end{array}
\right)  $, which proves the assertation.
\end{proof}

\textbf{Theorem 2.1. }\textit{The linear operator }$A$\textit{ is symmetric in
}$H$\textit{.}

\begin{proof}
Let $F,$ $G\in D\left(  A\right)  $. By two partial integrations, we get%
\begin{align}
\left\langle AF,G\right\rangle  &  =\left\langle F,AG\right\rangle +\left(
W\left(  f,\overline{g};h_{1}-0\right)  -W\left(  f,\overline{g};-1\right)
\right)  +\delta_{1}^{2}\left(  W\left(  f,\overline{g};h_{2}-0\right)
-W\left(  f,\overline{g};h_{1}+0\right)  \right) \nonumber\\
&  +\delta_{1}^{2}\delta_{2}^{2}\left(  W\left(  f,\overline{g};h_{3}%
-0\right)  -W\left(  f,\overline{g};h_{2}+0\right)  \right)  +...+%
{\displaystyle\prod\limits_{i=1}^{m}}
\delta_{i}^{2}\left(  W\left(  f,\overline{g};1\right)  -W\left(
f,\overline{g};h_{m}+0\right)  \right) \nonumber\\
&  +\frac{%
{\displaystyle\prod\limits_{i=1}^{m}}
\delta_{i}^{2}}{\rho}\left(  R_{1}^{^{\prime}}\left(  f\right)  R_{1}\left(
\overline{g}\right)  -R_{1}\left(  f\right)  R_{1}^{^{\prime}}\left(
\overline{g}\right)  \right)  \tag{6}%
\end{align}
where%
\[
W\left(  f,\overline{g};x\right)  =f\left(  x\right)  \overline{g}^{^{\prime}%
}\left(  x\right)  -f^{^{\prime}}\left(  x\right)  \overline{g}\left(
x\right)
\]
denotes the Wronskians of the functions $f$ \ and $\overline{g}$. Since $f$
\ and $\overline{g}$ satisfy the boundary condition (2), it follows that%
\begin{equation}
W\left(  f,\overline{g};-1\right)  =0. \tag{7}%
\end{equation}
From the transmission conditions (4)-(5) we get%
\begin{equation}
\left\{
\begin{array}
[c]{c}%
W\left(  f,\overline{g};h_{1}-0\right)  =\delta_{1}^{2}W\left(  f,\overline
{g};h_{1}+0\right)  ,\\
W\left(  f,\overline{g};h_{2}-0\right)  =\delta_{2}^{2}W\left(  f,\overline
{g};h_{2}+0\right)  ,\\
\vdots\\
W\left(  f,\overline{g};h_{m}-0\right)  =\delta_{m}^{2}W\left(  f,\overline
{g};h_{m}+0\right)  .
\end{array}
\right.  \tag{8}%
\end{equation}
Furthermore,%
\begin{align}
&  R_{1}^{^{\prime}}\left(  f\right)  R_{1}\left(  \overline{g}\right)
-R_{1}\left(  f\right)  R_{1}^{^{\prime}}\left(  \overline{g}\right)
\nonumber\\
&  =\left(  \beta_{1}^{^{\prime}}f\left(  1\right)  -\beta_{2}^{^{\prime}%
}f^{^{\prime}}\left(  1\right)  \right)  \left(  \beta_{1}\overline{g}\left(
1\right)  -\beta_{2}\overline{g}^{^{\prime}}\left(  1\right)  \right)
-\left(  \beta_{1}f\left(  1\right)  -\beta_{2}f^{^{\prime}}\left(  1\right)
\right)  \left(  \beta_{1}^{^{\prime}}\overline{g}\left(  1\right)  -\beta
_{2}^{^{\prime}}\overline{g}\left(  1\right)  \right) \nonumber\\
&  =\left(  \beta_{2}\beta_{1}^{^{\prime}}-\beta_{2}^{^{\prime}}\beta
_{1}\right)  f^{^{\prime}}\left(  1\right)  \overline{g}\left(  1\right)
+\left(  \beta_{2}^{^{\prime}}\beta_{1}-\beta_{2}\beta_{1}^{^{\prime}}\right)
f\left(  1\right)  \overline{g}^{^{\prime}}\left(  1\right) \nonumber\\
&  =\rho\left(  f^{^{\prime}}\left(  1\right)  \overline{g}\left(  1\right)
-f\left(  1\right)  \overline{g}^{^{\prime}}\left(  1\right)  \right)  =-\rho
W\left(  f,g;1\right)  . \tag{9}%
\end{align}
Finally, substituting (7)-(9) in (6) then we get
\begin{equation}
\left\langle AF,G\right\rangle =\left\langle F,AG\right\rangle \text{ }\left(
F,G\in D\left(  A\right)  \right)  \tag{10}%
\end{equation}

\end{proof}

Now we can write the following theorem with the helps of Theorem 2.1,
Naimark's Patching Lemma [15] and using the similar way in [6]

\textbf{Theorem 2.2. }\textit{The linear operator }$A$\textit{ is self-adjoint
in }$H$\textit{.}

\textbf{Corollary 2.1.} \textit{All eigenvalues of the problem (1)-(5) are
real.}

We can now assume that all eigenfunctions are real-valued.

\textbf{Corollary 2.2. }\textit{If }$\lambda_{1}$\textit{\ and }$\lambda_{2}%
$\textit{ are two different eigenvalues of the problem (1)-(5), then the
corresponding eigenfunctions }$u_{1}$\textit{ and }$u_{2}$\textit{ of this
problem satisfy the following equality}:%
\begin{align*}
&
{\displaystyle\int\limits_{-1}^{h_{1}}}
u_{1}\left(  x\right)  u_{2}\left(  x\right)  dx+\delta_{1}^{2}%
{\displaystyle\int\limits_{h_{1}}^{h_{2}}}
u_{1}\left(  x\right)  u_{2}\left(  x\right)  dx+\delta_{1}^{2}\delta_{2}^{2}%
{\displaystyle\int\limits_{h_{2}}^{h_{3}}}
u_{1}\left(  x\right)  u_{2}\left(  x\right)  dx+\\
&  ...+%
{\displaystyle\prod\limits_{i=1}^{m}}
\delta_{i}^{2}%
{\displaystyle\int\limits_{h_{m}}^{1}}
u_{1}\left(  x\right)  u_{2}\left(  x\right)  dx=-\frac{%
{\displaystyle\prod\limits_{i=1}^{m}}
\delta_{i}^{2}}{\rho}R_{1}^{^{\prime}}\left(  u_{1}\right)  R_{1}^{^{\prime}%
}\left(  u_{2}\right)  .
\end{align*}
\textit{In fact this formula means the orthogonality of eigenfunctions }%
$u_{1}$\textit{ and }$u_{2}$\textit{ in the Hilbert space }$H.$

We need the following lemma, which can be proved by the same technique as in
$\left[  4\right]  .$

\textbf{Lemma 2.3.} \textit{Let the real-valued function }$q\left(  x\right)
$\textit{ be continuous in }$\left[  -1,1\right]  $\textit{ and }$f\left(
\lambda\right)  ,g\left(  \lambda\right)  $\textit{ are given entire
functions. Then for any }$\lambda\in%
\mathbb{C}
$\textit{ the equation }%
\[
-u^{\prime\prime}+q\left(  x\right)  u=\lambda u,\text{ \ \ }x\in\left[
-1,1\right]
\]
\textit{has a unique solution }$u=u\left(  x,\lambda\right)  $\textit{
satisfies the initial conditions }%
\[
u\left(  -1\right)  =f\left(  \lambda\right)  ,u^{^{\prime}}\left(  -1\right)
=g\left(  \lambda\right)  \ \left(  \text{or }u\left(  1\right)  =f\left(
\lambda\right)  ,u^{^{\prime}}\left(  1\right)  =g\left(  \lambda\right)
\right)  .
\]
\textit{For each fixed }$x\in\left[  -1,1\right]  ,$\textit{ }$u\left(
x,\lambda\right)  $\textit{ is an entire function of }$\lambda.$

We shall define two solutions%
\[
\phi_{\lambda}\left(  x\right)  =\left\{
\begin{array}
[c]{c}%
\phi_{1\lambda}\left(  x\right)  ,\text{ }x\in\left[  -1,h_{1}\right)  ,\\
\phi_{2\lambda}\left(  x\right)  ,\text{ }x\in\left(  h_{1},{h}_{2}\right)
,\\
\vdots\\
\phi_{\left(  m+1\right)  \lambda}\left(  x\right)  ,\text{ }x\in\left(
h_{m},{1}\right]  ,
\end{array}
\right.  \text{ and }\chi_{\lambda}\left(  x\right)  =\left\{
\begin{array}
[c]{c}%
\chi_{1\lambda}\left(  x\right)  ,\text{ }x\in\left[  -1,h_{1}\right)  ,\\
\chi_{2\lambda}\left(  x\right)  ,\text{ }x\in\left(  h_{1},{h}_{2}\right)
,\\
\vdots\\
\chi_{\left(  m+1\right)  \lambda}\left(  x\right)  ,\text{ }x\in\left(
h_{m},{1}\right]  ,
\end{array}
\right.
\]
of the Eq. (1) as follows: Let $\phi_{1\lambda}\left(  x\right)  :=\phi
_{1}\left(  x,\lambda\right)  $ be the solution of Eq. (1) on $\left[
-1,h_{1}\right]  $, which satisfies the initial conditions%
\begin{equation}
u\left(  -1\right)  =\alpha_{2},\text{ }u^{\prime}\left(  -1\right)
=-\alpha_{1}. \tag{11}%
\end{equation}
By virtue of Lemma 2.1, after defining this solution, we may define the
solution $\phi_{2}\left(  x,\lambda\right)  :=\phi_{2\lambda}\left(  x\right)
$ of Eq. (1) on $\left[  h_{1},h_{2}\right]  $ by means of the solution
$\phi_{1}\left(  x,\lambda\right)  $ by the initial conditions%
\begin{equation}
u\left(  h_{1}\right)  =\delta_{1}^{-1}\phi_{1}\left(  h_{1},\lambda\right)
,\text{ }u^{\prime}\left(  h_{1}\right)  =\delta_{1}^{-1}\phi_{1}^{\prime
}\left(  h_{1},\lambda\right)  . \tag{12}%
\end{equation}
After defining this solution, we may define the solution $\phi_{3}\left(
x,\lambda\right)  :=\phi_{3\lambda}\left(  x\right)  $ of Eq. (1) on $\left[
h_{2},h_{3}\right]  $ by means of the solution $\phi_{2}\left(  x,\lambda
\right)  $ by the initial conditions%
\begin{equation}
u\left(  h_{2}\right)  =\delta_{2}^{-1}\phi_{2}\left(  h_{2},\lambda\right)
,\text{ }u^{\prime}\left(  h_{2}\right)  =\delta_{2}^{-1}\phi_{2}^{\prime
}\left(  h_{2},\lambda\right)  . \tag{13}%
\end{equation}
Continuing in this manner, we may define the solution $\phi_{\left(
m+1\right)  }\left(  x,\lambda\right)  :=\phi_{\left(  m+1\right)  \lambda
}\left(  x\right)  $ of Eq. (1) on $\left[  h_{m},1\right]  $ by means of the
solution $\phi_{m}\left(  x,\lambda\right)  $ by the initial conditions%
\begin{equation}
u\left(  h_{m}\right)  =\delta_{m}^{-1}\phi_{m}\left(  h_{m},\lambda\right)
,\text{ }u^{\prime}\left(  h_{m}\right)  =\delta_{m}^{-1}\phi_{m}^{\prime
}\left(  h_{m},\lambda\right)  . \tag{14}%
\end{equation}
Therefore, $\phi\left(  x,\lambda\right)  $ satisfies the Eq. (1) on $\left[
-1,h_{1}\right)  \cup\left(  h_{1},h_{2}\right)  \cup...\cup\left(
h_{m},1\right]  $, the boundary condition (3), and the transmission conditions (4)-(5).

Analogically, first we define the solution $\chi_{\left(  m+1\right)  \lambda
}\left(  x\right)  :=\chi_{\left(  m+1\right)  }\left(  x,\lambda\right)  $ on
$\left[  h_{m},1\right]  $ by the initial conditions%
\begin{equation}
u\left(  1\right)  =\beta_{2}^{\prime}\lambda+\beta_{2},\text{ }u^{\prime
}\left(  1\right)  =\beta_{1}^{\prime}\lambda+\beta_{1}. \tag{15}%
\end{equation}
Again, after defining this solution, we may define the solution $\chi
_{m\lambda}\left(  x\right)  :=\chi_{m}\left(  x,\lambda\right)  $ of the Eq.
(1) on $\left[  h_{m-1},h_{m}\right]  $ by the initial conditions%
\begin{equation}
u\left(  h_{m}\right)  =\delta_{m}\chi_{m+1}\left(  h_{m},\lambda\right)
,\text{ }u^{\prime}\left(  h_{m}\right)  =\delta_{m}\chi_{m+1}^{\prime}\left(
h_{m},\lambda\right)  . \tag{16}%
\end{equation}
Continuing in this manner, we may define the solution $\chi_{1\lambda}\left(
x\right)  :=\chi_{1}\left(  x,\lambda\right)  $ of the Eq. (1) on $\left[
-1,h_{1}\right]  $ by the initial conditions%
\begin{equation}
u\left(  h_{1}\right)  =\delta_{1}\chi_{2}\left(  h_{1},\lambda\right)
,\text{ }u^{^{\prime}}\left(  h_{1}\right)  =\delta_{1}\chi_{2}^{^{\prime}%
}\left(  h_{1},\lambda\right)  . \tag{17}%
\end{equation}
Therefore, $\chi\left(  x,\lambda\right)  $ satisfies the Eq. (1) on $\left[
-1,h_{1}\right)  \cup\left(  h_{1},h_{2}\right)  \cup...\cup\left(
h_{m},1\right]  $, the boundary condition (3), and the transmission conditions (4)-(5).

It is obvious that the Wronskians%
\[
\omega_{i}\left(  \lambda\right)  :=W_{\lambda}\left(  \phi_{i},\chi
_{i};x\right)  :=\phi_{i}\left(  x,\lambda\right)  \chi_{i}^{^{\prime}}\left(
x,\lambda\right)  -\phi_{i}^{^{\prime}}\left(  x,\lambda\right)  \chi
_{i}\left(  x,\lambda\right)  ,\text{ \ }x\in\Omega_{i}\text{ }\left(
i=\overline{1,m+1}\right)
\]
are independent of $x\in\Omega_{i}$ and entire functions.

\textbf{Lemma 2.4.} \textit{For each }$\lambda\in%
\mathbb{C}
$\textit{, }$\omega_{1}\left(  \lambda\right)  =\delta_{1}^{2}\omega
_{2}\left(  \lambda\right)  =\delta_{1}^{2}\delta_{2}^{2}\omega_{3}\left(
\lambda\right)  =...=\left(
{\displaystyle\prod\limits_{i=1}^{m}}
\delta_{i}^{2}\right)  \omega_{m+1}\left(  \lambda\right)  .$

\begin{proof}
By the means of (12), (13), (14), (16) and (17), the short calculation gives%
\begin{align*}
W_{\lambda}\left(  \phi_{1},\chi_{1};h_{1}\right)   &  =\delta_{1}%
^{2}W_{\lambda}\left(  \phi_{2},\chi_{2};h_{1}\right)  =\delta_{1}%
^{2}W_{\lambda}\left(  \phi_{3},\chi_{3};h_{2}\right)  =\delta_{1}^{2}%
\delta_{2}^{2}W_{\lambda}\left(  \phi_{3},\chi_{3};h_{2}\right) \\
&  =\delta_{1}^{2}\delta_{2}^{2}W_{\lambda}\left(  \phi_{4},\chi_{4}%
;h_{3}\right)  =...=\left(
{\displaystyle\prod\limits_{i=1}^{m}}
\delta_{i}^{2}\right)  W_{\lambda}\left(  \phi_{m+1},\chi_{m+1};h_{m}\right)
\end{align*}
so $\omega_{1}\left(  \lambda\right)  =\delta_{1}^{2}\omega_{2}\left(
\lambda\right)  =\delta_{1}^{2}\delta_{2}^{2}\omega_{3}\left(  \lambda\right)
=...=\left(
{\displaystyle\prod\limits_{i=1}^{m}}
\delta_{i}^{2}\right)  \omega_{m+1}\left(  \lambda\right)  .$
\end{proof}

Now we may introduce the characteristic function%
\[
\omega\left(  \lambda\right)  :=\omega_{1}\left(  \lambda\right)  =\delta
_{1}^{2}\omega_{2}\left(  \lambda\right)  =\delta_{1}^{2}\delta_{2}^{2}%
\omega_{3}\left(  \lambda\right)  =...=\left(
{\displaystyle\prod\limits_{i=1}^{m}}
\delta_{i}^{2}\right)  \omega_{m+1}\left(  \lambda\right)  .
\]

\textbf{Theorem 2.3.} \textit{The eigenvalues of the problem (1)-(5) are the
zeros of the function }$\omega\left(  \lambda\right)  .$

\begin{proof}
Let $\omega\left(  \lambda_{0}\right)  =0$. Then $W_{\lambda_{0}}\left(
\phi_{1},\chi_{1};x\right)  =0$ and therefore the functions $\phi
_{1\lambda_{0}}\left(  x\right)  $ and $\chi_{1\lambda_{0}}\left(  x\right)  $
are linearly dependent, i.e.%
\[
\chi_{1\lambda_{0}}\left(  x\right)  =k_{1}\phi_{1\lambda_{0}}\left(
x\right)  ,\text{ \ }x\in\left[  -1,h_{1}\right]
\]
for some $k_{1}\neq0$. From this, it follows that $\chi\left(  x,\lambda
_{0}\right)  $ satisfies also the first boundary condition (2), so
$\chi\left(  x,\lambda_{0}\right)  $ is an eigenfunction of the problem
$\left(  1\right)  -\left(  5\right)  $ corresponding to this eigenvalue
$\lambda_{0}$.

Now we let $u_{0}\left(  x\right)  $ be any eigenfunction corresponding to
eigenvalue $\lambda_{0}$, but $\omega\left(  \lambda_{0}\right)  \neq0$. Then
the functions $\phi_{1},\chi_{1},$ $\phi_{2},\chi_{2},...,\phi_{m+1}%
,\chi_{m+1}$ would be linearly independent on $\left[  -1,h_{1}\right]  ,$
$\left[  h_{1},h_{2}\right]  $ and $\left[  h_{m},1\right]  $ respectively.
Therefore $u_{0}\left(  x\right)  $ may be represented in the following form
\[
u_{0}\left(  x\right)  =\left\{
\begin{array}
[c]{c}%
c_{1}\phi_{1}\left(  x,\lambda_{0}\right)  +c_{2}\chi_{1}\left(  x,\lambda
_{0}\right)  ,\text{ \ }x\in\left[  -1,h_{1}\right)  ,\\
c_{3}\phi_{2}\left(  x,\lambda_{0}\right)  +c_{4}\chi_{2}\left(  x,\lambda
_{0}\right)  ,\text{ \ }x\in\left(  h_{1},h_{2}\right)  ,\\
\vdots\\
c_{2m+1}\phi_{m+1}\left(  x,\lambda_{0}\right)  +c_{2m+2}\chi_{m+1}\left(
x,\lambda_{0}\right)  ,\text{ \ }x\in\left(  h_{m},1\right]  .
\end{array}
\right.
\]
where at least one of the constants $c_{1},$ $c_{2},...,c_{2m+2}$ is not zero.
Considering the equations%
\begin{equation}
\tau_{\upsilon}\left(  u_{0}\left(  x\right)  \right)  =0,\text{ \ }%
\upsilon=\overline{1,2m+2} \tag{18}%
\end{equation}
as the homogenous system of linear equations of the variables $c_{1},$
$c_{2},$ $c_{2n+2}$ and taking (12), (13), (14), (16) and (17) into account,
it follows that the determinant of this system is%
\[
\left\vert
\begin{array}
[c]{cccccccc}%
0 & \omega_{1}\left(  \lambda_{0}\right)  & 0 & \cdots & \cdots & \cdots & 0 &
0\\
\phi_{1\lambda_{0}}\left(  h_{1}\right)  & \chi_{1\lambda_{0}}\left(
h_{1}\right)  & -\delta_{1}\phi_{2\lambda_{0}}\left(  h_{1}\right)  &
-\delta_{1}\phi_{2\lambda_{0}}\left(  h_{1}\right)  & \cdots & \cdots & 0 &
0\\
\phi_{1\lambda_{0}}^{\prime}\left(  h_{1}\right)  & \chi_{1\lambda_{0}%
}^{\prime}\left(  h_{1}\right)  & -\delta_{1}\phi_{2\lambda_{0}}^{\prime
}\left(  h_{1}\right)  & -\delta_{1}\phi_{2\lambda_{0}}^{\prime}\left(
h_{1}\right)  & \cdots & \cdots & \cdots & \vdots\\
0 & \ddots & \ddots & \ddots & \ddots & \ddots & \cdots & \vdots\\
\vdots & \vdots & \ddots & \ddots & \ddots & \ddots & \cdots & 0\\
\vdots & \vdots & \ddots & \ddots & \ddots & \ddots & -\delta_{m}\phi_{\left(
m+1\right)  \lambda_{0}}\left(  h_{m}\right)  & -\delta_{m}\chi_{\left(
m+1\right)  \lambda_{0}}\left(  h_{m}\right) \\
0 & 0 & \cdots & \cdots & \cdots & 0 & -\delta_{m}\phi_{\left(  m+1\right)
\lambda_{0}}^{\prime}\left(  h_{m}\right)  & -\delta_{m}\chi_{\left(
m+1\right)  \lambda_{0}}^{\prime}\left(  h_{m}\right) \\
0 & 0 & \cdots & \cdots & \cdots & 0 & \omega_{m+1}\left(  \lambda_{0}\right)
& 0
\end{array}
\right\vert
\]%
\[
=-\left(
{\displaystyle\prod\limits_{i=1}^{m}}
\delta_{i}^{2}\omega_{i}\left(  \lambda_{0}\right)  \right)  \omega_{m+1}%
^{m}\left(  \lambda_{0}\right)  \neq0.
\]
Therefore, the system (18) has only the trivial solution $c_{i}=0$ $\left(
i=\overline{1,2m+2}\right)  $. Thus we get a contradiction, which completes
the proof.
\end{proof}

\textbf{Lemma 2.5.} \textit{If }$\lambda=\lambda_{0}$\textit{ is an
eigenvalue, then }$\phi\left(  x,\lambda_{0}\right)  $\textit{ and }%
$\chi\left(  x,\lambda_{0}\right)  $\textit{ are linearly dependent.}

\begin{proof}
Let $\lambda=\lambda_{0}$ be an eigenvalue. Then by virtue of Theorem 2.2
\[
W\left(  \phi_{i\lambda_{0}},\chi_{i\lambda_{0}};x\right)  =\omega_{i}\left(
\lambda_{0}\right)  =0
\]
and hence%
\begin{equation}
\chi_{i\lambda_{0}}\left(  x\right)  =k_{i}\phi_{i\lambda_{0}}\left(
x\right)  \text{ \ \ }\left(  i=\overline{1,m+1}\right)  \tag{19}%
\end{equation}
for some $k_{1}\neq0,$ $k_{2}\neq0,...,k_{m+1}\neq0,$ We must show that
$k_{1}=k_{2}=...=k_{m+1}$. Suppose, if possible, that $k_{m}\neq k_{m+1}$.

Taking into account the definitions of the solutions $\phi_{i}\left(
x,\lambda\right)  $ and $\chi_{i}\left(  x,\lambda\right)  $ from the
equalities (19), we have%
\begin{align*}
\tau_{2m+1}\left(  \chi_{\lambda_{0}}\right)   &  =\chi_{\lambda_{0}}\left(
h_{m}-0\right)  -\delta_{m}\chi_{\lambda_{0}}\left(  h_{m}+0\right)
=\chi_{m\lambda_{0}}\left(  h_{m}\right)  -\delta_{m}\chi_{\left(  m+1\right)
\lambda_{0}}\left(  h_{m}\right) \\
&  =k_{m}\phi_{m}\left(  h_{m}\right)  -\delta_{m}k_{m+1}\phi_{m+1}\left(
h_{m}\right)  =k_{m}\delta_{m}\phi_{m+1}\left(  h_{m}\right)  -\delta
_{m}k_{m+1}\phi_{m+1}\left(  h_{m}\right) \\
&  =\delta_{m}\left(  k_{m}-k_{m+1}\right)  \phi_{m+1}\left(  h_{m}\right)
=0.
\end{align*}
since $\tau_{2m+1}\left(  \chi_{\lambda_{0}}\right)  =0$ and $\delta
_{m}\left(  k_{m}-k_{m+1}\right)  \neq0,$ it follows that%
\begin{equation}
\phi_{\left(  m+1\right)  \lambda_{0}}\left(  h_{m}\right)  =0. \tag{20}%
\end{equation}
By the same procedure from $\tau_{2n+2}\left(  \chi_{\lambda_{0}}\right)  =0$
we can derive that%
\begin{equation}
\phi_{\left(  m+1\right)  \lambda_{0}}^{^{\prime}}\left(  h_{m}\right)  =0.
\tag{21}%
\end{equation}
From the fact that $\phi_{\left(  m+1\right)  \lambda_{0}}\left(  x\right)  $
is a solution of the differential equation (1) on $\left[  h_{m},1\right]  $
and satisfies the initial conditions (20) and (21), it follows that
$\phi_{\left(  m+1\right)  \lambda_{0}}\left(  x\right)  =0$ identically on
$\left[  h_{m},1\right]  $ because of the well-known existence and uniqueness
theorem for the initial value problems of the ordinary linear differential
equations. Making use of (14), (19) and (20), we may also derive that%
\begin{equation}
\phi_{m\lambda_{0}}\left(  h_{m}\right)  =\phi_{m\lambda_{0}}^{^{\prime}%
}\left(  h_{m}\right)  =0 \tag{22}%
\end{equation}
Continuing in this matter, we may also find that%
\begin{equation}
\left\{
\begin{array}
[c]{c}%
\phi_{\left(  m-1\right)  \lambda_{0}}\left(  h_{m-1}\right)  =\phi_{\left(
m-1\right)  \lambda_{0}}^{^{\prime}}\left(  h_{m-1}\right)  =0.\\
\vdots\\
\phi_{1\lambda_{0}}\left(  h_{1}\right)  =\phi_{1\lambda_{0}}^{^{\prime}%
}\left(  h_{1}\right)  =0.
\end{array}
\right.  \tag{23}%
\end{equation}
identically on $\left[  h_{m-1},h_{m}\right]  ,...,\left[  -1,h_{1}\right]  $
respectively. Hence $\phi\left(  x,\lambda_{0}\right)  =0$ identically on
$\left[  -1,h_{1}\right)  \cup\left(  h_{1},h_{2}\right)  \cup...\cup\left(
h_{m},1\right]  $. But this contradicts with (11). Hence $k_{m}=k_{m+1}.$
Analogically we can prove that $k_{m-1}=k_{m},$...,$k_{2}=k_{3}$ and
$k_{1}=k_{2}.$
\end{proof}

\textbf{Corollary 2.3.} \textit{If }$\lambda=\lambda_{0}$\textit{\ is an
eigenvalue, then both }$\phi\left(  x,\lambda_{0}\right)  $\textit{ and }%
$\chi\left(  x,\lambda_{0}\right)  $\textit{ are eigenfunctions corresponding
to this eigenvalue.}

\textbf{Lemma 2.6.} \textit{All eigenvalues }$\lambda_{n}$\textit{ are simple
zeros of }$\omega(\lambda).$

\begin{proof}
Using the Lagrange's formula (cf. [25], p. 6-7), it can be shown that%
\begin{equation}
\left(  \lambda-\lambda_{n}\right)  \left(
{\displaystyle\int\limits_{-1}^{h_{1}}}
\phi_{\lambda}\left(  x\right)  \phi_{\lambda_{n}}\left(  x\right)
dx+\delta_{1}^{2}%
{\displaystyle\int\limits_{h_{1}}^{h_{2}}}
\phi_{\lambda}\left(  x\right)  \phi_{\lambda_{n}}\left(  x\right)
dx+\delta_{1}^{2}\delta_{2}^{2}%
{\displaystyle\int\limits_{h_{2}}^{h_{3}}}
\phi_{\lambda}\left(  x\right)  \phi_{\lambda_{n}}\left(  x\right)  dx+\right.
\nonumber
\end{equation}%
\begin{equation}
\left.  ...+\left(
{\displaystyle\prod\limits_{i=1}^{m}}
\delta_{i}^{2}\right)
{\displaystyle\int\limits_{h_{m}}^{1}}
\phi_{\lambda}\left(  x\right)  \phi_{\lambda_{n}}\left(  x\right)  dx\right)
=\left(
{\displaystyle\prod\limits_{i=1}^{m}}
\delta_{i}^{2}\right)  W\left(  \phi_{\lambda},\phi_{\lambda_{n}};1\right)
\tag{24}%
\end{equation}
for any $\lambda$. Recall that%
\[
\chi_{\lambda_{n}}\left(  x\right)  =k_{n}\phi_{\lambda_{n}}\left(  x\right)
,\text{ \ }x\in\left[  -1,h_{1}\right)  \cup\left(  h_{1},h_{2}\right)
\cup...\cup\left(  h_{m},1\right]
\]
for some $k_{n}\neq0,$ $n=1,2,...$. Using this equality for the right side of
(24), we have%
\begin{align*}
W\left(  \phi_{\lambda},\phi_{\lambda_{n}};1\right)   &  =\frac{1}{k_{n}%
}W\left(  \phi_{\lambda},\chi_{\lambda_{n}};1\right)  =\frac{1}{k_{n}}\left(
\lambda_{n}R_{1}^{^{\prime}}\left(  \phi_{\lambda}\right)  +R_{1}\left(
\phi_{\lambda}\right)  \right) \\
&  =\frac{1}{k_{n}}\left[  \omega\left(  \lambda\right)  -\left(
\lambda-\lambda_{n}\right)  R_{1}^{^{\prime}}\left(  \phi_{\lambda}\right)
\right] \\
&  =\left(  \lambda-\lambda_{n}\right)  \frac{1}{k_{n}}\left[  \frac
{\omega\left(  \lambda\right)  }{\lambda-\lambda_{n}}-R_{1}^{^{\prime}}\left(
\phi_{\lambda}\right)  \right]  .
\end{align*}
Substituting this formula in (24) and letting $\lambda\rightarrow\lambda_{n}$,
we get%
\begin{align}
&
{\displaystyle\int\limits_{-1}^{h_{1}}}
\left(  \phi_{\lambda_{n}}\left(  x\right)  \right)  ^{2}dx+\delta_{1}^{2}%
{\displaystyle\int\limits_{h_{1}}^{h_{2}}}
\left(  \phi_{\lambda_{n}}\left(  x\right)  \right)  ^{2}dx+\delta_{1}%
^{2}\delta_{2}^{2}%
{\displaystyle\int\limits_{h_{2}}^{h_{3}}}
\left(  \phi_{\lambda_{n}}\left(  x\right)  \right)  ^{2}dx+...+%
{\displaystyle\prod\limits_{i=1}^{m}}
\delta_{i}^{2}%
{\displaystyle\int\limits_{h_{m}}^{1}}
\left(  \phi_{\lambda_{n}}\left(  x\right)  \right)  ^{2}dx\nonumber\\
&  =\frac{%
{\displaystyle\prod\limits_{i=1}^{m}}
\delta_{i}^{2}}{k_{n}}\left(  \omega^{^{\prime}}\left(  \lambda_{n}\right)
-R_{1}^{^{\prime}}\left(  \phi_{\lambda_{n}}\right)  \right)  . \tag{25}%
\end{align}
Now putting%
\[
R_{1}^{^{\prime}}\left(  \phi_{\lambda_{n}}\right)  =\frac{1}{k_{n}}%
R_{1}^{^{\prime}}\left(  \chi_{\lambda_{n}}\right)  =\frac{\rho}{k_{n}}%
\]
in (25) we get $\omega^{^{\prime}}\left(  \lambda_{n}\right)  \neq0.$
\end{proof}

\textbf{Definition 2.1. }\textit{The geometric multiplicity of an eigenvalue
}$\lambda$ \textit{of the problem (1)-(5) is the dimension of its eigenspace,
i.e. the number of its linearly independent eigenfunctions.}

\textbf{Theorem 2.4. }\textit{All eigenvalues of the problem (1)-(5) are
geometrically simple.}

\begin{proof}
If $f$ and $g$ are two eigenfunctions for an eigenvalue $\lambda_{0}$ of
(1)-(5) then (2) implies that $f(-1)=cg(-1)$ and $f^{\prime}(-1)=cg^{\prime
}(-1)$ for some constant $c\in%
\mathbb{C}
$. By the uniqueness theorem for solutions of ordinary differential equation
and the transmission conditions (4)-(5), we have that $f=cg$ on $\left[
-1,h_{1}\right]  ,$ $\left[  h_{1},h_{2}\right]  $ and $\left[  h_{m}%
,1\right]  $. Thus the geometric multiplicity of $\lambda_{0}$ is one.
\end{proof}

\section{Asymptotic approximate formulas of $\omega\left(  \lambda\right)  $
for four distinct cases}

We start by proving some lemmas.

\textbf{Lemma 3.1.} \textit{Let }$\phi\left(  x,\lambda\right)  $\textit{ be
the solutions of Eq. (1) defined in Section 2, and let }$\lambda=s^{2}%
$\textit{.Then the following integral equations hold for }$k=0,1:$%
\begin{equation}
\left\{
\begin{array}
[c]{c}%
\phi_{1\lambda}^{^{\left(  k\right)  }}\left(  x\right)  =\alpha_{2}\left(
\cos s\left(  x+1\right)  \right)  ^{\left(  k\right)  }-\alpha_{1}\frac{1}%
{s}\left(  \sin s\left(  x+1\right)  \right)  ^{\left(  k\right)  }\\
+\frac{1}{s}%
{\displaystyle\int\limits_{-1}^{x}}
\left(  \sin s\left(  x-y\right)  \right)  ^{\left(  k\right)  }q\left(
y\right)  \phi_{1\lambda}\left(  y\right)  dy,\\
\phi_{2\lambda}^{^{\left(  k\right)  }}\left(  x\right)  =\frac{1}{\delta_{1}%
}\phi_{1\lambda}\left(  h_{1}\right)  \left(  \cos s\left(  x-h_{1}\right)
\right)  ^{\left(  k\right)  }+\frac{1}{s}\frac{1}{\delta_{1}}\phi_{1\lambda
}^{^{\prime}}\left(  h_{1}\right)  \left(  \sin s\left(  x-h_{1}\right)
\right)  ^{\left(  k\right)  }\\
+\frac{1}{s}%
{\displaystyle\int\limits_{h_{1}}^{x}}
\left(  \sin s\left(  x-y\right)  \right)  ^{\left(  k\right)  }q\left(
y\right)  \phi_{2\lambda}\left(  y\right)  dy,\\
\vdots\\
\phi_{\left(  m+1\right)  \lambda}^{^{\left(  k\right)  }}\left(  x\right)
=\frac{1}{\delta_{m}}\phi_{m\lambda}\left(  h_{m}\right)  \left(  \cos
s\left(  x-h_{m}\right)  \right)  ^{\left(  k\right)  }+\frac{1}{s}\frac
{1}{\delta_{m}}\phi_{m\lambda}^{^{\prime}}\left(  h_{m}\right)  \left(  \sin
s\left(  x-h_{m}\right)  \right)  ^{\left(  k\right)  }\\
+\frac{1}{s}%
{\displaystyle\int\limits_{h_{m}}^{x}}
\left(  \sin s\left(  x-y\right)  \right)  ^{\left(  k\right)  }q\left(
y\right)  \phi_{\left(  m+1\right)  \lambda}\left(  y\right)  dy,
\end{array}
\right.  \tag{26}%
\end{equation}
\textit{where }$\left(  \cdot\right)  ^{\left(  k\right)  }=\frac{d^{k}%
}{dx^{k}}\left(  \cdot\right)  $\textit{.}

\begin{proof}
It is enough to substitute $s^{2}\phi_{1\lambda}\left(  y\right)
+\phi_{1\lambda}^{^{\prime\prime}}\left(  y\right)  ,$ $s^{2}\phi_{2\lambda
}\left(  y\right)  +\phi_{2\lambda}^{^{\prime\prime}}\left(  y\right)
,...,s^{2}\phi_{\left(  m+1\right)  \lambda}\left(  y\right)  +\phi_{\left(
m+1\right)  \lambda}^{^{\prime\prime}}\left(  y\right)  $ instead of $q\left(
y\right)  \phi_{1\lambda}\left(  y\right)  ,$ $q\left(  y\right)
\phi_{2\lambda}\left(  y\right)  ,$ $q\left(  y\right)  \phi_{\left(
m+1\right)  \lambda}\left(  y\right)  $ in the integral terms of the (26),
respectively, and integrate by parts twice.
\end{proof}

\textbf{Lemma 3.2.} \textit{Let }$\lambda=s^{2}$\textit{, }$\operatorname{Im}%
s=t.$\textit{ Then the functions }$\phi_{i\lambda}\left(  x\right)  $\textit{
have the following asymptotic formulas for }$\left\vert \lambda\right\vert
\rightarrow\infty,$\textit{ which hold uniformly for }$x\in\Omega_{i}$
$\left(  \text{\textit{for} }i=\overline{1,m+1}\text{ and }k=0,1.\right)  :$%
\begin{align}
\phi_{1\lambda}^{^{\left(  k\right)  }}\left(  x\right)   &  =\alpha
_{2}\left(  \cos s\left(  x+1\right)  \right)  ^{\left(  k\right)  }+O\left(
\left\vert s\right\vert ^{k-1}e^{\left\vert t\right\vert \left(  x+1\right)
}\right)  ,\nonumber\\
\phi_{2\lambda}^{^{\left(  k\right)  }}\left(  x\right)   &  =\frac{\alpha
_{2}}{\delta_{1}}\left(  \cos s\left(  x+1\right)  \right)  ^{\left(
k\right)  }+O\left(  \left\vert s\right\vert ^{k-1}e^{\left\vert t\right\vert
\left(  x+1\right)  }\right)  ,\nonumber\\
&  \vdots\nonumber\\
\phi_{\left(  m+1\right)  \lambda}^{^{\left(  k\right)  }}\left(  x\right)
&  =\frac{\alpha_{2}}{%
{\displaystyle\prod\limits_{i=1}^{m}}
\delta_{i}}\left(  \cos s\left(  x+1\right)  \right)  ^{\left(  k\right)
}+O\left(  \left\vert s\right\vert ^{k-1}e^{\left\vert t\right\vert \left(
x+1\right)  }\right)  \tag{27}%
\end{align}
\textit{if} $\alpha_{2}\neq0,$%
\begin{align}
\phi_{1\lambda}^{^{\left(  k\right)  }}\left(  x\right)   &  =-\frac
{\alpha_{1}}{s}\left(  \sin s\left(  x+1\right)  \right)  ^{\left(  k\right)
}+O\left(  \left\vert s\right\vert ^{k-2}e^{\left\vert t\right\vert \left(
x+1\right)  }\right)  ,\nonumber\\
\phi_{2\lambda}^{^{\left(  k\right)  }}\left(  x\right)   &  =-\frac
{\alpha_{1}}{s\delta_{1}}\left(  \sin s\left(  x+1\right)  \right)  ^{\left(
k\right)  }+O\left(  \left\vert s\right\vert ^{k-2}e^{\left\vert t\right\vert
\left(  x+1\right)  }\right)  ,\nonumber\\
&  \vdots\nonumber\\
\phi_{\left(  m+1\right)  \lambda}^{^{\left(  k\right)  }}\left(  x\right)
&  =-\frac{\alpha_{1}}{s%
{\displaystyle\prod\limits_{i=1}^{m}}
\delta_{i}}\left(  \sin s\left(  x+1\right)  \right)  ^{\left(  k\right)
}+O\left(  \left\vert s\right\vert ^{k-2}e^{\left\vert t\right\vert \left(
x+1\right)  }\right)  \tag{28}%
\end{align}
\textit{if} $\alpha_{2}=0.$

\begin{proof}
Since the proof of the formulae for $\phi_{1\lambda}\left(  x\right)  $ is
identical to Titchmarsh's proof to similar results for $\phi_{\lambda}\left(
x\right)  $ (see [26], Lemma 1.7 p. 9-10), we may formulate them without
proving them here.

Since the proof of the formulae for $\phi_{2\lambda}\left(  x\right)  $ and
$\phi_{3\lambda}\left(  x\right)  $ are identical to Kadakal's and Mukhtarov's
proof to similar results for $\phi_{\lambda}\left(  x\right)  $ (see [11],
Lemma 3.2 p. 1373-1375), we may formulate them without proving them here. But
the similar formulae for $\phi_{4\lambda}\left(  x\right)  ,...,\phi_{\left(
m+1\right)  \lambda}(x)$ need individual consideration, since the last
solutions are defined by the initial conditions of these special nonstandart
forms. We shall only prove the formula $\left(  27\right)  $ for $k=0$ and
$m=3$.

Let $\alpha_{2}\neq0$. Then according to $(27)$ for $m=2$%
\begin{align*}
\phi_{3\lambda}\left(  h_{3}\right)   &  =\alpha_{2}\left\{  \frac{1}%
{\delta_{2}}\cos s\left(  h_{3}-h_{2}\right)  \left[  \frac{1}{\delta_{1}}\cos
s\left(  h_{2}-h_{1}\right)  \cos s\left(  h_{1}+1\right)  \right.  \right. \\
&  \left.  -\frac{1}{\delta_{1}}\sin s\left(  h_{2}-h_{1}\right)  \sin
s\left(  h_{1}+1\right)  \right]  -\frac{1}{\delta_{2}}\sin s\left(
h_{3}-h_{2}\right)  \left[  \frac{1}{\delta_{1}}\sin s\left(  h_{2}%
-h_{1}\right)  \cos s\left(  h_{1}+1\right)  \right. \\
&  \left.  \left.  +\frac{1}{\delta_{1}}\cos s\left(  h_{2}-h_{1}\right)  \sin
s\left(  h_{1}+1\right)  \right]  \right\}  +O\left(  \left\vert s\right\vert
^{-1}e^{\left\vert t\right\vert \left[  \left(  h_{3}-h_{2}\right)  +\left(
h_{2}-h_{1}\right)  +\left(  h_{1}+1\right)  \right]  }\right)
\end{align*}
and%
\begin{align*}
\phi_{3\lambda}^{^{\prime}}\left(  h_{3}\right)   &  =\alpha_{2}\left\{
-\frac{s}{\delta_{2}}\sin s\left(  h_{3}-h_{2}\right)  \left[  \frac{1}%
{\delta_{1}}\cos s\left(  h_{2}-h_{1}\right)  \cos s\left(  h_{1}+1\right)
\right.  \right. \\
&  \left.  -\frac{1}{\delta_{1}}\sin s\left(  h_{2}-h_{1}\right)  \sin
s\left(  h_{1}+1\right)  \right]  -\frac{s}{\delta_{2}}\cos s\left(
h_{3}-h_{2}\right)  \left[  \frac{1}{\delta_{1}}\sin s\left(  h_{2}%
-h_{1}\right)  \cos s\left(  h_{1}+1\right)  \right. \\
&  \left.  \left.  +\frac{1}{\delta_{1}}\cos s\left(  h_{2}-h_{1}\right)  \sin
s\left(  h_{1}+1\right)  \right]  \right\}  +O\left(  e^{\left\vert
t\right\vert \left[  \left(  h_{3}-h_{2}\right)  +\left(  h_{2}-h_{1}\right)
+\left(  h_{1}+1\right)  \right]  }\right)
\end{align*}
Substituting these asymptotic expressions into $\left(  26\right)  $, we get%
\begin{align}
\phi_{4\lambda}\left(  x\right)   &  =\alpha_{2}\left\{  \frac{1}{\delta_{3}%
}\left(  \cos s\left(  x-h_{3}\right)  \right)  \left\{  \frac{1}{\delta_{2}%
}\cos s\left(  h_{3}-h_{2}\right)  \right.  \left[  \frac{1}{\delta_{1}}\cos
s\left(  h_{2}-h_{1}\right)  \cos s\left(  h_{1}+1\right)  \right.  \right.
\nonumber\\
&  \left.  -\frac{1}{\delta_{1}}\sin s\left(  h_{2}-h_{1}\right)  \sin
s\left(  h_{1}+1\right)  \right]  -\frac{1}{\delta_{2}}\sin s\left(
h_{3}-h_{2}\right)  \left[  \frac{1}{\delta_{1}}\sin s\left(  h_{2}%
-h_{1}\right)  \cos s\left(  h_{1}+1\right)  \right. \nonumber\\
&  \left.  \left.  +\frac{1}{\delta_{1}}\cos s\left(  h_{2}-h_{1}\right)  \sin
s\left(  h_{1}+1\right)  \right]  \right\}  -\frac{1}{\delta_{3}}\sin s\left(
x-h_{3}\right)  \left\{  \frac{1}{\delta_{2}}\sin s\left(  h_{3}-h_{2}\right)
\right. \nonumber\\
&  \times\left[  \frac{1}{\delta_{1}}\cos s\left(  h_{2}-h_{1}\right)  \cos
s\left(  h_{1}+1\right)  -\frac{1}{\delta_{1}}\sin s\left(  h_{2}%
-h_{1}\right)  \sin s\left(  h_{1}+1\right)  \right]  +\frac{1}{\delta_{2}%
}\nonumber\\
&  \left.  \left.  \times\cos s\left(  h_{3}-h_{2}\right)  \left[  \frac
{1}{\delta_{1}}\sin s\left(  h_{2}-h_{1}\right)  \cos s\left(  h_{1}+1\right)
+\frac{1}{\delta_{1}}\cos s\left(  h_{2}-h_{1}\right)  \sin s\left(
h_{1}+1\right)  \right]  \right\}  \right\} \nonumber\\
&  +\frac{1}{s}%
{\displaystyle\int\limits_{h_{3}}^{x}}
\sin s\left(  x-y\right)  q\left(  y\right)  \phi_{4\lambda}\left(  y\right)
dy+O\left(  \left\vert s\right\vert ^{-1}e^{\left\vert t\right\vert \left[
\left(  x-h_{3}\right)  +\left(  h_{3}-h_{2}\right)  +\left(  h_{2}%
-h_{1}\right)  +\left(  h_{1}+1\right)  \right]  }\right)  \tag{29}%
\end{align}
Multiplying through by $e^{-\left\vert t\right\vert \left[  \left(
x-h_{3}\right)  +\left(  h_{3}-h_{2}\right)  +\left(  h_{2}-h_{1}\right)
+\left(  h_{1}+1\right)  \right]  }$, and denoting%
\[
F_{4\lambda}\left(  x\right)  :=e^{-\left\vert t\right\vert \left[  \left(
x-h_{3}\right)  +\left(  h_{3}-h_{2}\right)  +\left(  h_{2}-h_{1}\right)
+\left(  h_{1}+1\right)  \right]  }\phi_{4\lambda}\left(  x\right)
\]
we have%
\begin{align*}
F_{4\lambda}\left(  x\right)   &  :=\alpha_{2}e^{-\left\vert t\right\vert
\left[  \left(  x-h_{3}\right)  +\left(  h_{3}-h_{2}\right)  +\left(
h_{2}-h_{1}\right)  +\left(  h_{1}+1\right)  \right]  }\left\{  \frac
{1}{\delta_{3}}\cos s\left(  x-h_{3}\right)  \left\{  \frac{1}{\delta_{2}}\cos
s\left(  h_{3}-h_{2}\right)  \right.  \right. \\
&  \times\left[  \frac{1}{\delta_{1}}\cos s\left(  h_{2}-h_{1}\right)  \cos
s\left(  h_{1}+1\right)  -\frac{1}{\delta_{1}}\sin s\left(  h_{2}%
-h_{1}\right)  \sin s\left(  h_{1}+1\right)  \right]  -\frac{1}{\delta_{2}%
}\sin s\left(  h_{3}-h_{2}\right) \\
&  \left.  \times\left[  \frac{1}{\delta_{1}}\sin s\left(  h_{2}-h_{1}\right)
\cos s\left(  h_{1}+1\right)  +\frac{1}{\delta_{1}}\cos s\left(  h_{2}%
-h_{1}\right)  \sin s\left(  h_{1}+1\right)  \right]  \right\}  -\frac
{1}{\delta_{3}}\sin s\left(  x-h_{3}\right) \\
&  \times\left\{  \frac{1}{\delta_{2}}\sin s\left(  h_{3}-h_{2}\right)
\right.  \left[  \frac{1}{\delta_{1}}\cos s\left(  h_{2}-h_{1}\right)  \cos
s\left(  h_{1}+1\right)  -\frac{1}{\delta_{1}}\sin s\left(  h_{2}%
-h_{1}\right)  \sin s\left(  h_{1}+1\right)  \right] \\
&  +\frac{1}{\delta_{2}}\cos s\left(  h_{3}-h_{2}\right)  \left[  \frac
{1}{\delta_{1}}\sin s\left(  h_{2}-h_{1}\right)  \cos s\left(  h_{1}+1\right)
+\frac{1}{\delta_{1}}\cos s\left(  h_{2}-h_{1}\right)  \sin s\left(
h_{1}+1\right)  \right] \\
&  +\frac{1}{s}%
{\displaystyle\int\limits_{h_{3}}^{x}}
\sin s\left(  x-y\right)  q\left(  y\right)  e^{-\left\vert t\right\vert
\left[  \left(  x-h_{3}\right)  +\left(  h_{3}-h_{2}\right)  +\left(
h_{2}-h_{1}\right)  +\left(  h_{1}+1\right)  \right]  }F_{4\lambda}\left(
y\right)  dy+O\left(  \left\vert s\right\vert ^{-1}\right)  .
\end{align*}
Denoting $M:=\max_{x\in\left[  h_{3},1\right]  }\left\vert F_{4\lambda}\left(
x\right)  \right\vert $ from the last formula, it follows that%
\[
M\left(  \lambda\right)  \leq\left\vert \alpha_{2}\right\vert \frac
{2}{\left\vert \delta_{1}\right\vert }\frac{2}{\left\vert \delta
_{2}\right\vert }\frac{2}{\left\vert \delta_{3}\right\vert }+\frac{M\left(
\lambda\right)  }{\left\vert s\right\vert }%
{\displaystyle\int\limits_{h_{3}}^{1}}
q\left(  y\right)  dy+\frac{M_{0}}{\left\vert s\right\vert }%
\]
for some $M_{0}>0$. From this, it follows that $M\left(  \lambda\right)
=O\left(  1\right)  $ as $\lambda\rightarrow\infty$, so%
\[
\phi_{4\lambda}\left(  x\right)  =O\left(  e^{\left\vert t\right\vert \left[
\left(  x-h_{3}\right)  +\left(  h_{3}-h_{2}\right)  +\left(  h_{2}%
-h_{1}\right)  +\left(  h_{1}+1\right)  \right]  }\right)  .
\]
Substituting this back into the integral on the right side of (29) yields (27)
for $k=0$ and $m=3$. The other cases may be considered analogically.
\end{proof}

\textbf{Theorem 3.1.} \textit{Let }$\lambda=s^{2}$\textit{, }%
$t=\operatorname{Im}s$\textit{. Then the characteristic function }%
$\omega\left(  \lambda\right)  $\textit{ has the following asymptotic
formulas}$:$

\textit{Case 1}$:$\textit{ If }$\beta_{2}^{\prime}\neq0$\textit{, }$\alpha
_{2}\neq0$\textit{, then}%
\begin{equation}
\omega\left(  \lambda\right)  =\beta_{2}^{\prime}\alpha_{2}s^{3}\left(
{\displaystyle\prod\limits_{i=1}^{m}}
\delta_{i}^{2}\right)  \sin2s+O\left(  \left\vert s\right\vert ^{2}%
e^{2\left\vert t\right\vert }\right)  . \tag{30}%
\end{equation}
\textit{Case 2}$:$\textit{ If }$\beta_{2}^{\prime}\neq0$\textit{, }$\alpha
_{2}=0$\textit{, then}%
\begin{equation}
\omega\left(  \lambda\right)  =-\beta_{2}^{\prime}\alpha_{1}s^{2}\left(
{\displaystyle\prod\limits_{i=1}^{m}}
\delta_{i}^{2}\right)  \cos2s+O\left(  \left\vert s\right\vert ^{2}%
e^{2\left\vert t\right\vert }\right)  . \tag{31}%
\end{equation}
\textit{Case 3}$:$\textit{ If }$\beta_{2}^{\prime}=0$\textit{, }$\alpha
_{2}\neq0$\textit{, then}%
\begin{equation}
\omega\left(  \lambda\right)  =\beta_{1}^{\prime}\alpha_{2}s^{2}\left(
{\displaystyle\prod\limits_{i=1}^{m}}
\delta_{i}^{2}\right)  \cos2s+O\left(  \left\vert s\right\vert ^{2}%
e^{2\left\vert t\right\vert }\right)  . \tag{32}%
\end{equation}
\textit{Case 4}$:$\textit{ If }$\beta_{2}^{\prime}=0$\textit{, }$\alpha_{2}%
=0$\textit{, then}%
\begin{equation}
\omega\left(  \lambda\right)  =-\beta_{1}^{\prime}\alpha_{1}s\left(
{\displaystyle\prod\limits_{i=1}^{m}}
\delta_{i}^{2}\right)  \sin2s+O\left(  \left\vert s\right\vert ^{2}%
e^{2\left\vert t\right\vert }\right)  . \tag{33}%
\end{equation}

\begin{proof}
The proof is completed by substituting (27) and (28) into the representation%
\begin{align}
\omega\left(  \lambda\right)   &  =\left(
{\displaystyle\prod\limits_{i=1}^{m}}
\delta_{i}^{2}\right)  \omega_{m+1}\left(  \lambda\right)  =\left(
{\displaystyle\prod\limits_{i=1}^{m}}
\delta_{i}^{2}\right)  \left[  \phi_{\left(  m+1\right)  \lambda}\left(
1\right)  \chi_{\left(  m+1\right)  \lambda}\left(  1\right)  -\phi_{\left(
m+1\right)  \lambda}^{^{\prime}}\left(  1\right)  \chi_{\left(  m+1\right)
\lambda}^{^{\prime}}\left(  1\right)  \right]  =\nonumber\\
&  =\left(
{\displaystyle\prod\limits_{i=1}^{m}}
\delta_{i}^{2}\right)  \left[  \left(  \lambda\beta_{1}^{^{\prime}}+\beta
_{1}\right)  \phi_{\left(  m+1\right)  \lambda}\left(  1\right)  -\left(
\lambda\beta_{2}^{^{\prime}}+\beta_{2}\right)  \phi_{\left(  m+1\right)
\lambda}^{^{\prime}}\left(  1\right)  \right] \nonumber\\
&  =\left(
{\displaystyle\prod\limits_{i=1}^{m}}
\delta_{i}^{2}\right)  \left[  \lambda\left(  \beta_{1}^{\prime}\phi_{\left(
m+1\right)  \lambda}\left(  1\right)  -\beta_{2}^{\prime}\phi_{\left(
m+1\right)  \lambda}^{\prime}\left(  1\right)  \right)  +\left(  \beta_{1}%
\phi_{\left(  m+1\right)  \lambda}\left(  1\right)  -\beta_{2}\phi_{\left(
m+1\right)  \lambda}^{^{\prime}}\left(  1\right)  \right)  \right] \nonumber\\
&  =-\lambda\left(
{\displaystyle\prod\limits_{i=1}^{m}}
\delta_{i}^{2}\right)  \beta_{2}^{^{\prime}}\phi_{\left(  m+1\right)  \lambda
}^{^{\prime}}\left(  1\right)  +\lambda\left(
{\displaystyle\prod\limits_{i=1}^{m}}
\delta_{i}^{2}\right)  \beta_{1}^{^{\prime}}\phi_{\left(  m+1\right)  \lambda
}\left(  1\right)  +\left(
{\displaystyle\prod\limits_{i=1}^{m}}
\delta_{i}^{2}\right)  -\beta_{2}\phi_{\left(  m+1\right)  \lambda}^{^{\prime
}}\left(  1\right) \nonumber\\
&  +\left(
{\displaystyle\prod\limits_{i=1}^{m}}
\delta_{i}^{2}\right)  \beta_{1}\phi_{\left(  m+1\right)  \lambda}\left(
1\right)  . \tag{34}%
\end{align}

\end{proof}

\textbf{Corollary 3.1.} \textit{The eigenvalues of the problem (1)-(5) are
bounded below.}

\begin{proof}
Putting $s=it$ $\left(  t>0\right)  $ in the above formulas, it follows that
$\omega\left(  -t^{2}\right)  \rightarrow\infty$ as $t\rightarrow\infty$.
Therefore, $\omega\left(  \lambda\right)  \neq0$ for $\lambda$ negative and
sufficiently large.
\end{proof}

\section{Asymptotic formulas for eigenvalues and eigenfunctions}

Now we can obtain the asymptotic approximation formulae for the eigenvalues of
the considered problem (1)-(5).

Since the eigenvalues coincide with the zeros of the entire function
$\omega_{m+1}\left(  \lambda\right)  $, it follows that they have no finite
limit. Moreover, we know from Corollaries 2.1 and 3.1 that all eigenvalues are
real and bounded below. Hence, we may renumber them as $\lambda_{0}\leq
\lambda_{1}\leq\lambda_{2}\leq...$, listed according to their multiplicity.

\textbf{Theorem 4.1.} \textit{The eigenvalues }$\lambda_{n}=s_{n}^{2}%
$\textit{, }$n=0,1,2,...$\textit{ of the problem (1)-(5) have the following
asymptotic formulae for }$n\rightarrow\infty:$

\textit{Case 1}$:$\textit{ If }$\beta_{2}^{^{\prime}}\neq0,\alpha_{2}\neq
0,$\textit{ then}%
\begin{equation}
s_{n}=\frac{\pi\left(  n-1\right)  }{2}+O\left(  \frac{1}{n}\right)  .
\tag{35}%
\end{equation}

\textit{Case 2}$:$\textit{ If }$\beta_{2}^{^{\prime}}\neq0,\alpha_{2}%
=0,$\textit{ then}%
\begin{equation}
s_{n}=\frac{\pi\left(  n-\frac{1}{2}\right)  }{2}+O\left(  \frac{1}{n}\right)
. \tag{36}%
\end{equation}

\textit{Case 3}$:$\textit{ If }$\beta_{2}^{^{\prime}}=0,\alpha_{2}\neq
0,$\textit{ then}%
\begin{equation}
s_{n}=\frac{\pi\left(  n-\frac{1}{2}\right)  }{2}+O\left(  \frac{1}{n}\right)
. \tag{37}%
\end{equation}

\textit{Case 4}$:$\textit{ If }$\beta_{2}^{^{\prime}}=0,\alpha_{2}=0,$\textit{
then}%
\begin{equation}
s_{n}=\frac{\pi n}{2}+O\left(  \frac{1}{n}\right)  . \tag{38}%
\end{equation}

\begin{proof}
We shall only consider the first case. The other cases may be considered
similarly. Denoting $\omega_{1}\left(  s\right)  $ and $\omega_{2}\left(
s\right)  $ the first and $O$-term of the right of $\left(  42\right)  $
repectively, we shall apply the well-known Rouch\'{e}'s theorem, which asserts
that if $f\left(  s\right)  $ and $g\left(  s\right)  $ are analytic inside
and on a closed contour $C$, and $\left\vert g\left(  s\right)  \right\vert
<\left\vert f\left(  s\right)  \right\vert $ on $C$, then $f\left(  s\right)
$ and $f\left(  s\right)  +g\left(  s\right)  $ have the same number zeros
inside $C$, provided that each zero is counted according to their
multiplicity. It is readily shown that $\left\vert \overline{\omega}%
_{1}\left(  s\right)  \right\vert >\left\vert \overline{\omega}_{2}\left(
s\right)  \right\vert $ on the contours%
\[
C_{n}:=\left\{  s\in%
\mathbb{C}
\left\vert \left\vert s\right\vert \right.  =\frac{\left(  n+\frac{1}%
{2}\right)  \pi}{2}\right\}
\]
for sufficiently large $n$.

Let $\lambda_{0}\leq\lambda_{1}\leq\lambda_{2}\leq...$ be zeros of
$\omega\left(  \lambda\right)  $ and $\lambda_{n}=s_{n}^{2}$. Since inside the
contour $C_{n}$, $\overline{\omega}_{1}\left(  s\right)  $ has zeros at points
$s=0$ and $s=\frac{k\pi}{4},$ $k=\pm1,$ $\pm2,...,$ $\pm n$.%
\begin{equation}
s_{n}=\frac{\left(  n-1\right)  \pi}{2}+\delta_{n} \tag{39}%
\end{equation}
where $\delta_{n}=O\left(  1\right)  $ for sufficiently large $n$. By
substituting this in (30), we derive that $\delta_{n}=O\left(  \frac{1}%
{n}\right)  $, which completes the proof.
\end{proof}

The next approximation for the eigenvalues may be obtained by the following
procedure. For this, we shall suppose that $q\left(  y\right)  $ is of bounded
variation in $\left[  -1,1\right]  $.

Firstly we consider the case $\beta_{2}^{^{\prime}}\neq0$ and $\alpha_{2}%
\neq0$. Putting $x=h_{1},x=h_{2},...,x=h_{m}$ in $\left(  26\right)  $ and
then substituting in the expression of $\phi_{\left(  m+1\right)  \lambda
}^{^{\prime}}$, we get that%
\begin{align*}
\phi_{\left(  m+1\right)  \lambda}^{^{\prime}}\left(  1\right)   &
=s\frac{\alpha_{2}}{%
{\displaystyle\prod\limits_{i=1}^{m}}
\delta_{i}}\sin2s-\frac{\alpha_{1}}{%
{\displaystyle\prod\limits_{i=1}^{m}}
\delta_{i}}\cos2s+\frac{1}{\left(
{\displaystyle\prod\limits_{i=1}^{m}}
\delta_{i}\right)  }%
{\displaystyle\int\limits_{-1}^{h_{1}}}
\cos\left(  s\left(  1-y\right)  \right)  q\left(  y\right)  \phi_{1\lambda
}\left(  y\right)  dy\\
&  +\frac{1}{\left(
{\displaystyle\prod\limits_{i=2}^{m}}
\delta_{i}\right)  }%
{\displaystyle\int\limits_{h_{1}}^{h_{2}}}
\cos\left(  s\left(  1-y\right)  \right)  q\left(  y\right)  \phi_{2\lambda
}\left(  y\right)  dy+\frac{1}{\left(
{\displaystyle\prod\limits_{i=3}^{m}}
\delta_{i}\right)  }%
{\displaystyle\int\limits_{h_{2}}^{h_{3}}}
\cos\left(  s\left(  1-y\right)  \right)  q\left(  y\right)  \phi_{3\lambda
}\left(  y\right)  dy\\
&  +...+\frac{1}{\delta_{m}}%
{\displaystyle\int\limits_{h_{m-1}}^{h_{m}}}
\cos\left(  s\left(  1-y\right)  \right)  q\left(  y\right)  \phi_{m\lambda
}\left(  y\right)  dy.+%
{\displaystyle\int\limits_{h_{m}}^{1}}
\cos\left(  s\left(  1-y\right)  \right)  q\left(  y\right)  \phi_{\left(
m+1\right)  \lambda}\left(  y\right)  dy.
\end{align*}

Substituting (27) into the right side of the last integral equality then gives%
\begin{align*}
\phi_{\left(  m+1\right)  \lambda}^{^{\prime}}\left(  1\right)   &
=-\frac{s\alpha_{2}}{%
{\displaystyle\prod\limits_{i=1}^{m}}
\delta_{i}}\sin2s-\frac{\alpha_{1}}{%
{\displaystyle\prod\limits_{i=1}^{m}}
\delta_{i}}\cos2s+\frac{\alpha_{2}}{\left(
{\displaystyle\prod\limits_{i=1}^{m}}
\delta_{i}\right)  }%
{\displaystyle\int\limits_{-1}^{h_{1}}}
\cos\left(  s\left(  1-y\right)  \right)  \cos\left(  s\left(  1+y\right)
\right)  q\left(  y\right)  dy\\
&  +\frac{\alpha_{2}}{\left(
{\displaystyle\prod\limits_{i=1}^{m}}
\delta_{i}\right)  }%
{\displaystyle\int\limits_{h_{1}}^{h_{2}}}
\cos\left(  s\left(  1-y\right)  \right)  \cos\left(  s\left(  1+y\right)
\right)  q\left(  y\right)  dy+\\
&  ...+\frac{\alpha_{2}}{\left(
{\displaystyle\prod\limits_{i=1}^{m}}
\delta_{i}\right)  }%
{\displaystyle\int\limits_{h_{m}}^{1}}
\cos\left(  s\left(  1-y\right)  \right)  \cos\left(  s\left(  1+y\right)
\right)  q\left(  y\right)  dy+O\left(  \left\vert s\right\vert ^{-1}%
e^{2\left\vert t\right\vert }\right)  .
\end{align*}
On the other hand, from (27), it follows that%
\[
\phi_{\left(  m+1\right)  \lambda}\left(  1\right)  =\frac{\alpha_{2}}{%
{\displaystyle\prod\limits_{i=1}^{m}}
\delta_{i}}\cos2s+O\left(  \left\vert s\right\vert ^{-1}e^{2\left\vert
t\right\vert }\right)  .
\]
Putting these formulas into (34), we have%
\begin{align*}
\omega\left(  \lambda\right)   &  =\frac{s^{3}\beta_{2}^{^{\prime}}\alpha_{2}%
}{%
{\displaystyle\prod\limits_{i=1}^{m}}
\delta_{i}}\sin2s+s^{2}\left[  \left(  \frac{\beta_{1}^{^{\prime}}\alpha
_{2}+\beta_{2}^{^{\prime}}\alpha_{1}}{%
{\displaystyle\prod\limits_{i=1}^{m}}
\delta_{i}}\right)  \cos2s\right. \\
&  -\frac{\beta_{2}^{^{\prime}}}{\left(
{\displaystyle\prod\limits_{i=1}^{m}}
\delta_{i}\right)  }%
{\displaystyle\int\limits_{-1}^{h_{1}}}
\cos\left(  s\left(  1-y\right)  \right)  q\left(  y\right)  \phi_{1\lambda
}\left(  y\right)  dy-\frac{\beta_{2}^{^{\prime}}}{\left(
{\displaystyle\prod\limits_{i=2}^{m}}
\delta_{i}\right)  }%
{\displaystyle\int\limits_{h_{1}}^{h_{2}}}
\cos\left(  s\left(  1-y\right)  \right)  q\left(  y\right)  \phi_{2\lambda
}\left(  y\right)  dy-\\
&  \left.  ...-\frac{\beta_{2}^{^{\prime}}}{\delta_{m}}%
{\displaystyle\int\limits_{h_{m-1}}^{h_{m}}}
\cos\left(  s\left(  1-y\right)  \right)  q\left(  y\right)  \phi_{m\lambda
}\left(  y\right)  dy-\beta_{2}^{^{\prime}}%
{\displaystyle\int\limits_{h_{m}}^{1}}
\cos\left(  s\left(  1-y\right)  \right)  q\left(  y\right)  \phi_{\left(
m+1\right)  \lambda}\left(  y\right)  dy\right] \\
&  +O\left(  \left\vert s\right\vert ^{-1}e^{2\left\vert t\right\vert
}\right)  .
\end{align*}

Putting (39) in the last equality we find that%
\begin{align}
\sin\left(  2\delta_{n}\right)   &  =-\frac{\cos\left(  2\delta_{n}\right)
}{s_{n}}\left[  \frac{\beta_{1}^{^{\prime}}}{\beta_{2}^{^{\prime}}}%
+\frac{\alpha_{1}}{\alpha_{2}}-\frac{1}{2\left(
{\displaystyle\prod\limits_{i=1}^{m}}
\delta_{i}\right)  }%
{\displaystyle\int\limits_{-1}^{1}}
q\left(  y\right)  dy-\frac{1}{2\left(
{\displaystyle\prod\limits_{i=1}^{m}}
\delta_{i}\right)  }%
{\displaystyle\int\limits_{-1}^{1}}
\cos\left(  2s_{n}y\right)  q\left(  y\right)  dy\right] \nonumber\\
&  +O\left(  \left\vert s_{n}\right\vert ^{-2}\right)  \tag{40}%
\end{align}
Recalling that $q\left(  y\right)  $ is of bounded variation in $\left[
-1,1\right]  $, and applying the well-known Riemann-Lebesque Lemma (see [27],
p. 48, Theorem 4.12) to the second integral on the right in (40), this term is
$O\left(  \frac{1}{n}\right)  $. As a result, from (40) it follows that%
\[
\delta_{n}=-\frac{1}{\pi\left(  n-1\right)  }\left[  \frac{\beta_{1}^{\prime}%
}{\beta_{2}^{\prime}}+\frac{\alpha_{1}}{\alpha_{2}}-\frac{1}{2\left(
{\displaystyle\prod\limits_{i=1}^{m}}
\delta_{i}\right)  }%
{\displaystyle\int\limits_{-1}^{1}}
q\left(  y\right)  dy\right]  +O\left(  \frac{1}{n^{2}}\right)  .
\]
Substituting in (40), we have%
\[
s_{n}=\frac{\pi\left(  n-1\right)  }{2}-\frac{1}{\pi\left(  n-1\right)
}\left[  \frac{\beta_{1}^{\prime}}{\beta_{2}^{\prime}}+\frac{\alpha_{1}%
}{\alpha_{2}}-\frac{1}{2\left(
{\displaystyle\prod\limits_{i=1}^{m}}
\delta_{i}\right)  }%
{\displaystyle\int\limits_{-1}^{1}}
q\left(  y\right)  dy\right]  +O\left(  \frac{1}{n^{2}}\right)  .
\]

Similar formulas in the other cases are as follows:

In case 2:%
\[
s_{n}=\frac{\pi\left(  n-\frac{1}{2}\right)  }{2}-\frac{1}{\pi\left(
n-\frac{1}{2}\right)  }\left[  \frac{\beta_{1}^{\prime}}{\beta_{2}^{\prime}%
}+\frac{1}{2\left(
{\displaystyle\prod\limits_{i=1}^{m}}
\delta_{i}\right)  }%
{\displaystyle\int\limits_{-1}^{1}}
q\left(  y\right)  dy\right]  +O\left(  \frac{1}{n^{2}}\right)  .
\]

In case 3:%
\[
s_{n}=\frac{\pi\left(  n-\frac{1}{2}\right)  }{2}+\frac{1}{\pi\left(
n-\frac{1}{2}\right)  }\left[  \frac{\beta_{2}}{\beta_{1}^{\prime}}%
-\frac{\alpha_{1}}{\alpha_{2}}+\frac{1}{2\left(
{\displaystyle\prod\limits_{i=1}^{m}}
\delta_{i}\right)  }%
{\displaystyle\int\limits_{-1}^{1}}
q\left(  y\right)  dy\right]  +O\left(  \frac{1}{n^{2}}\right)  .
\]

In case 4:%
\[
s_{n}=\frac{\pi n}{2}+\frac{1}{\pi n}\left[  \frac{\beta_{2}}{\beta
_{1}^{\prime}}+\frac{1}{2\left(
{\displaystyle\prod\limits_{i=1}^{m}}
\delta_{i}\right)  }%
{\displaystyle\int\limits_{-1}^{1}}
q\left(  y\right)  dy\right]  +O\left(  \frac{1}{n^{2}}\right)  .
\]
Recalling that $\phi\left(  x,\lambda_{n}\right)  $ is an eigenfunction
according to the eigenvalue $\lambda_{n}$ and by putting (35) into the (27) we
obtain that%
\[
\phi_{1\lambda_{n}}\left(  x\right)  =\alpha_{2}\cos\left(  \frac{\pi\left(
n-1\right)  \left(  x+1\right)  }{2}\right)  +O\left(  \frac{1}{n}\right)  ,
\]%
\begin{align*}
\phi_{2\lambda_{n}}\left(  x\right)   &  =\frac{\alpha_{2}}{\delta_{1}}%
\cos\left(  \frac{\pi\left(  n-1\right)  \left(  x+1\right)  }{2}\right)
+O\left(  \frac{1}{n}\right)  ,\\
&  \vdots
\end{align*}%
\[
\phi_{\left(  m+1\right)  \lambda_{n}}\left(  x\right)  =\frac{\alpha_{2}}{%
{\displaystyle\prod\limits_{i=1}^{m}}
\delta_{i}}\cos\left(  \frac{\pi\left(  n-1\right)  \left(  x+1\right)  }%
{2}\right)  +O\left(  \frac{1}{n}\right)
\]
in the first case. Consequently, if $\beta_{2}^{^{\prime}}\neq0$ and
$\alpha_{2}\neq0$, then the eigenfunction $\phi\left(  x,\lambda_{n}\right)  $
has the following asymptotic formulae%
\[
\phi\left(  x,\lambda_{n}\right)  =\left\{
\begin{array}
[c]{cc}%
\alpha_{2}\cos\left(  \frac{\pi\left(  n-1\right)  \left(  x+1\right)  }%
{2}\right)  +O\left(  \frac{1}{n}\right)  , & x\in\left[  -1,h_{1}\right) \\
\frac{\alpha_{2}}{\delta_{1}}\cos\left(  \frac{\pi\left(  n-1\right)  \left(
x+1\right)  }{2}\right)  +O\left(  \frac{1}{n}\right)  , & x\in\left(
h_{1},h_{2}\right) \\
\vdots & \vdots\\
\frac{\alpha_{2}}{%
{\displaystyle\prod\limits_{i=1}^{m}}
\delta_{i}}\cos\left(  \frac{\pi\left(  n-1\right)  \left(  x+1\right)  }%
{2}\right)  +O\left(  \frac{1}{n}\right)  , & x\in\left(  h_{m},1\right]
\end{array}
\right.
\]
which holds uniformly for $x\in\left[  -1,h_{1}\right)  \cup\left(
h_{1},h_{2}\right)  \cup...\cup\left(  h_{m},1\right]  $.

Similar formulas in the other cases are as follows:

In case 2%
\[
\phi\left(  x,\lambda_{n}\right)  =\left\{
\begin{array}
[c]{cc}%
-\frac{2\alpha_{1}}{\pi\left(  n-\frac{1}{2}\right)  }\sin\left(  \frac
{\pi\left(  n-\frac{1}{2}\right)  \left(  x+1\right)  }{2}\right)  +O\left(
\frac{1}{n^{2}}\right)  , & x\in\left[  -1,h_{1}\right) \\
-\frac{2\alpha_{1}}{\delta_{1}\pi\left(  n-\frac{1}{2}\right)  }\sin\left(
\frac{\pi\left(  n-\frac{1}{2}\right)  \left(  x+1\right)  }{2}\right)
+O\left(  \frac{1}{n^{2}}\right)  , & x\in\left(  h_{1},h_{2}\right) \\
\vdots & \vdots\\
-\frac{2\alpha_{1}}{\left(
{\displaystyle\prod\limits_{i=1}^{m}}
\delta_{i}\right)  \pi\left(  n-\frac{1}{2}\right)  }\sin\left(  \frac
{\pi\left(  n-\frac{1}{2}\right)  \left(  x+1\right)  }{2}\right)  +O\left(
\frac{1}{n^{2}}\right)  , & x\in\left(  h_{m},1\right]  .
\end{array}
\right.
\]

In case 3%
\[
\phi\left(  x,\lambda_{n}\right)  =\left\{
\begin{array}
[c]{cc}%
\alpha_{2}\cos\left(  \frac{\pi\left(  n-\frac{1}{2}\right)  \left(
x+1\right)  }{2}\right)  +O\left(  \frac{1}{n}\right)  , & x\in\left[
-1,h_{1}\right) \\
\frac{\alpha_{2}}{\delta_{1}}\cos\left(  \frac{\pi\left(  n-\frac{1}%
{2}\right)  \left(  x+1\right)  }{2}\right)  +O\left(  \frac{1}{n}\right)  , &
x\in\left(  h_{1},h_{2}\right) \\
\vdots & \vdots\\
\frac{\alpha_{2}}{%
{\displaystyle\prod\limits_{i=1}^{m}}
\delta_{i}}\cos\left(  \frac{\pi\left(  n-\frac{1}{2}\right)  \left(
x+1\right)  }{2}\right)  +O\left(  \frac{1}{n}\right)  , & x\in\left(
h_{m},1\right]  .
\end{array}
\right.
\]

In case 4%
\[
\phi\left(  x,\lambda_{n}\right)  =\left\{
\begin{array}
[c]{cc}%
-\frac{2\alpha_{1}}{\pi n}\sin\left(  \frac{\pi n\left(  x+1\right)  }%
{2}\right)  +O\left(  \frac{1}{n^{2}}\right)  , & x\in\left[  -1,h_{1}\right)
\\
-\frac{2\alpha_{1}}{\delta_{1}\pi n}\sin\left(  \frac{\pi n\left(  x+1\right)
}{2}\right)  +O\left(  \frac{1}{n^{2}}\right)  , & x\in\left(  h_{1}%
,h_{2}\right) \\
\vdots & \vdots\\
-\frac{2\alpha_{1}}{\left(
{\displaystyle\prod\limits_{i=1}^{m}}
\delta_{i}\right)  \pi n}\sin\left(  \frac{\pi n\left(  x+1\right)  }%
{2}\right)  +O\left(  \frac{1}{n^{2}}\right)  , & x\in\left(  h_{m},1\right]
.
\end{array}
\right.
\]
All these asymptotic formulas hold uniformly for $x\in\left[  -1,h_{1}\right)
\cup\left(  h_{1},h_{2}\right)  \cup...\cup\left(  h_{m},1\right]  $.

\section{Completeness of eigenfunctions}

\textbf{Theorem 5.1 }\textit{The operator }$A$\textit{ has only point
spectrum, i.e., }$\sigma\left(  A\right)  =\sigma_{\rho}\left(  A\right)  .$

\begin{proof}
It suffices to prove that if $\eta$ is not an eigenvalue of $A,$ then $\eta
\in\sigma\left(  A\right)  .$ Since $A$ is self-adjoint, we only consider a
real $\eta$. We investigate the equation $\left(  A-\eta\right)  Y=F\in H,$
where $F=\left(
\begin{array}
[c]{c}%
f\\
f_{1}%
\end{array}
\right)  .$

Let us consider the initial-value problem%
\begin{equation}
\left\{
\begin{array}
[c]{c}%
\tau y-\eta y=f,\text{ \ }x\in\left[  -1,h_{1}\right)  \cup\left(  h_{1}%
,h_{2}\right)  \cup...\cup\left(  h_{m},1\right]  ,\\
\alpha_{1}y\left(  -1\right)  +\alpha_{2}y^{\prime}\left(  -1\right)  =0,\\
y\left(  h_{i}-0\right)  -\delta_{i}y\left(  h_{i}+0\right)  =0,\\
y^{\prime}\left(  h_{i}-0\right)  -\delta_{i}y^{\prime}\left(  h_{i}+0\right)
=0.
\end{array}
\right.  \tag{41}%
\end{equation}
Let $u\left(  x\right)  $ be the solution of the equation $\tau u-\eta u=0$
satisfying%
\[
\left\{
\begin{array}
[c]{c}%
u\left(  -1\right)  =\alpha_{2},\text{ }u^{\prime}\left(  -1\right)
=-\alpha_{1},\\
u\left(  h_{i}-0\right)  -\delta_{i}u\left(  h_{i}+0\right)  =0,\\
u^{\prime}\left(  h_{i}-0\right)  -\delta_{i}u^{\prime}\left(  h_{i}+0\right)
=0.
\end{array}
\right.
\]
In fact,%
\[
u\left(  x\right)  =\left\{
\begin{array}
[c]{c}%
u_{1}\left(  x\right)  ,\text{ \ }x\in\left[  -1,h_{1}\right)  ,\\
u_{2}\left(  x\right)  ,\text{ \ }x\in\left(  h_{1},h_{2}\right)  ,\\
\vdots\\
u_{m+1}\left(  x\right)  ,\text{ \ }x\in\left(  h_{m},1\right]  ,
\end{array}
\right.
\]
where $u_{1}\left(  x\right)  $ is the unique solution of the initial-value
problem%
\[
\left\{
\begin{array}
[c]{c}%
-u^{\prime\prime}+q(x)u=\eta u,\text{ \ \ }x\in\left[  -1,h_{1}\right)  ,\\
u\left(  -1\right)  =\alpha_{2},\text{ }u^{\prime}\left(  -1\right)
=-\alpha_{1};
\end{array}
\right.
\]
$u_{2}\left(  x\right)  $ is the unique solution of the problem%
\[
\left\{
\begin{array}
[c]{c}%
-u^{\prime\prime}+q(x)u=\eta u,\text{ \ \ }x\in\left(  h_{1},h_{2}\right)  ,\\
u\left(  h_{1}-0\right)  =\delta_{1}u\left(  h_{1}+0\right)  ,\\
u^{\prime}\left(  h_{1}-0\right)  =\delta_{1}u^{\prime}\left(  h_{1}+0\right)
;
\end{array}
\right.
\]
$u_{3}\left(  x\right)  $ is the unique solution of the problem%
\[
\left\{
\begin{array}
[c]{c}%
-u^{\prime\prime}+q(x)u=\eta u,\text{ \ \ }x\in\left(  h_{2},h_{3}\right)  ,\\
u\left(  h_{2}-0\right)  =\delta_{2}u\left(  h_{2}+0\right)  ,\\
u^{\prime}\left(  h_{2}-0\right)  =\delta_{2}u^{\prime}\left(  h_{2}+0\right)
,
\end{array}
\right.
\]%
\[
\vdots
\]
and $u_{m+1}\left(  x\right)  $ is the unique solution of the problem%
\[
\left\{
\begin{array}
[c]{c}%
-u^{\prime\prime}+q(x)u=\eta u,\text{ \ \ }x\in\left(  h_{m},1\right]  ,\\
u\left(  h_{m}-0\right)  =\delta_{m}u\left(  h_{m}+0\right)  ,\\
u^{\prime}\left(  h_{m}-0\right)  =\delta_{m}u^{\prime}\left(  h_{m}+0\right)
.
\end{array}
\right.
\]
Let%
\[
w\left(  x\right)  =\left\{
\begin{array}
[c]{c}%
w_{1}\left(  x\right)  ,\text{ \ }x\in\left[  -1,h_{1}\right)  ,\\
w_{2}\left(  x\right)  ,\text{ \ }x\in\left(  h_{1},h_{2}\right)  ,\\
\vdots\\
w_{m+1}\left(  x\right)  ,\text{ \ }x\in\left(  h_{m},1\right]  ,
\end{array}
\right.
\]
be a solution of $\tau w-\eta w=f$ satisfying%
\[%
\begin{array}
[c]{c}%
\alpha_{1}w\left(  -1\right)  +\alpha_{2}w^{^{\prime}}\left(  -1\right)  =0,\\
w\left(  h_{i}-0\right)  =\delta_{i}w\left(  h_{i}+0\right)  ,\\
w^{^{\prime}}\left(  h_{i}-0\right)  =\delta_{i}w^{^{\prime}}\left(
h_{i}+0\right)  .
\end{array}
\]
Then, $\left(  41\right)  $ has the general solution%
\begin{equation}
y\left(  x\right)  =\left\{
\begin{array}
[c]{c}%
du_{1}+w_{1},\text{ \ }x\in\left[  -1,h_{1}\right)  ,\\
du_{2}+w_{2},\text{ \ }x\in\left(  h_{1},h_{2}\right)  ,\\
\vdots\\
du_{m+1}+w_{m+1},\text{ \ }x\in\left(  h_{m},1\right]  ,
\end{array}
\right.  \tag{42}%
\end{equation}
where $d\in%
\mathbb{C}
.$

Since $\eta$ is not an eigenvalue of the problem $\left(  1\right)  -\left(
5\right)  ,$ we have%
\begin{equation}
\eta\left[  \beta_{1}^{\prime}u_{m+1}\left(  1\right)  +\beta_{2}^{\prime
}u_{m+1}^{\prime}(1)\right]  +\left[  \beta_{1}u_{m+1}\left(  1\right)
+\beta_{2}u_{m+1}^{\prime}(1)\right]  \neq0. \tag{43}%
\end{equation}
The second component of $\left(  A-\eta\right)  Y=F$ involves the equation%
\[
-R_{1}\left(  y\right)  -\eta R_{1}^{\prime}\left(  y\right)  =f_{1},
\]
namely,%
\begin{equation}
\left[  -\beta_{1}y\left(  1\right)  -\beta_{2}y^{\prime}(1)\right]
-\eta\left[  \beta_{1}^{\prime}y\left(  1\right)  +\beta_{2}^{\prime}%
y^{\prime}(1)\right]  =f_{1}. \tag{44}%
\end{equation}
Substituting $\left(  42\right)  $ into $\left(  44\right)  $, we get%
\begin{align*}
&  \left(  -\beta_{2}u_{m+1}^{\prime}\left(  1\right)  -\beta_{1}%
u_{m+1}\left(  1\right)  -\eta\beta_{2}^{\prime}u_{m+1}^{\prime}\left(
1\right)  -\eta\beta_{1}^{\prime}u_{m+1}\left(  1\right)  \right)  d\\
&  =f_{1}+\beta_{1}w_{m+1}\left(  1\right)  +\beta_{2}w_{m+1}^{\prime}\left(
1\right)  +\eta\beta_{1}^{\prime}w_{m+1}\left(  1\right)  +\eta\beta
_{2}^{\prime}w_{m+1}^{\prime}\left(  1\right)
\end{align*}
In view of $\left(  43\right)  $, we know that $d$ is uniquely solvable.
Therefore, $y$ is uniquely determined.

The above arguments show that $\left(  A-\eta I\right)  ^{-1}$ is defined on
all of $H,$ where $I$ is identity matrix. We obtain that $\left(  A-\eta
I\right)  ^{-1}$ is bounded by Theorem 2.2 and the Closed Graph Theorem. Thus,
$\eta\in\sigma\left(  A\right)  .$ Therefore, $\sigma\left(  A\right)
=\sigma_{\rho}\left(  A\right)  .$
\end{proof}

The following lemma may be easily proved.

\textbf{Lemma 5.1 }\textit{The eigenvalues of the boundary value problem
}$\left(  1\right)  -\left(  5\right)  $\textit{ are bounded below, and they
are countably infinite and can cluster only at }$\infty.$

For every $\delta\in%
\mathbb{R}
\setminus\sigma_{\rho}\left(  A\right)  ,$ we have the following immediate conclusion.

\textbf{Lemma 5.2 }\textit{Let }$\lambda$\textit{ be an eigenvalue of
}$A-\delta I,$\textit{ and }$V$\textit{ a corresponding eigenfunction. Then,
}$\lambda^{-1}$\textit{ is an eigenvalue of }$\left(  A-\delta I\right)
^{-1},$\textit{ and }$V$\textit{ is a corresponding eigenfunction. The
converse is also true.}

On the other hand, if $\mu$ is an eigenvalue of $A$ and $U$ is a corresponding
eigenfunction, then $\mu-\delta$ is an eigenvalue of $A-\delta I,$ and $U$ is
a corresponding eigenfunction. The converse is also true. Accordingly, the
discussion about the completeness of the eigenfunctions of $A$ is equivalent
to considering the corresponding property of $\left(  A-\delta I\right)
^{-1}$.

By Lemma 1.1, Lemma 3.1 and Corollary 1.1, we suppose that $\left\{
\lambda_{n};\text{ }n\in%
\mathbb{N}
\right\}  $ is the real sequence of eigenvalues of $A,$ then $\left\{
\lambda_{n}-\delta;\text{ }n\in%
\mathbb{N}
\right\}  $ is the sequence of eigenvalues of $A-\delta I$. We may assume that%
\[
\left\vert \lambda_{1}-\delta\right\vert \leq\left\vert \lambda_{2}%
-\delta\right\vert \leq...\leq\left\vert \lambda_{n}-\delta\right\vert
\leq...\rightarrow\infty.
\]
Let $\left\{  \mu_{n};\text{ }n\in%
\mathbb{N}
\right\}  $ be the sequence of eigenvalues of $\left(  A-\delta I\right)
^{-1}.$ Then $\mu_{n}=\left(  \lambda_{n}-\delta\right)  ^{-1}$ and%
\[
\left\vert \mu_{1}\right\vert \geq\left\vert \mu_{2}\right\vert \geq
...\geq\left\vert \mu_{n}\right\vert \geq...\rightarrow0.
\]
Note that $0$ is not an eigenvalue of $\left(  A-\delta I\right)  ^{-1}.$

\textbf{Theorem 5.2 }\textit{The operator }$A$\textit{ has compact resolvents,
i.e, for each }$\delta\in%
\mathbb{R}
\setminus\sigma_{\rho}\left(  A\right)  ,$\textit{ }$\left(  A-\delta
I\right)  ^{-1}$\textit{ is compact on }$H.$

\begin{proof}
Let $\left\{  \mu_{1},\mu_{2},...\right\}  $ be the eigenvalues of $\left(
A-\delta I\right)  ^{-1},$ and $\left\{  P_{1},P_{2},...\right\}  $ the
orthogonal projections of finite rank onto the corresponding eigenspaces.
Since $\left\{  \mu_{1},\mu_{2},...\right\}  $ is a bounded sequence and all
$P_{n}\prime s$ are mutually orthogonal, we have $\sum_{n=1}^{\infty}\mu
_{n}P_{n}$ is strongly convergent to the bounded operator $\left(  A-\delta
I\right)  ^{-1},$ i.e., $\left(  A-\delta I\right)  ^{-1}=\sum_{n=1}^{\infty
}\mu_{n}P_{n}.$ Because for every $\alpha>0,$ the number of $\mu_{n}\prime s$
satisfying $\left\vert \mu_{n}\right\vert >\alpha$ is finite, and all
$P_{n}\prime s$ are of finite rank, we obtain that $\left(  A-\delta I\right)
^{-1}$ is compact.
\end{proof}

In terms of the above statements and the spectral theorem for compact
operators, we obtain the following theorem.

\textbf{Theorem 5.3 }\textit{The eigenfunctions of the problem }$\left(
1\right)  -\left(  5\right)  ,$\textit{ augmented to become eigenfunctions of
}$A,$\textit{ are complete in }$H$\textit{, i.e., if we let}%
\[
\left\{  \Phi_{n}=\left(
\begin{array}
[c]{c}%
\phi_{n}\left(  x\right) \\
R_{1}^{^{\prime}}\left(  \phi_{n}\right)
\end{array}
\right)  ;\text{ }n\in%
\mathbb{N}
\right\}
\]
\textit{ be a maximum set of orthonormal eigenfunctions of }$A,$\textit{ where
}$\left\{  \phi_{n}\left(  x\right)  ;\text{ }n\in%
\mathbb{N}
\right\}  $\textit{ are eigenfunctions of }$\left(  1\right)  -\left(
5\right)  ,$\textit{ then for all }$F\in H,$\textit{ }$F=\sum_{n=1}^{\infty
}\left\langle F,\Phi_{n}\right\rangle \Phi_{n}.$

\end{document}